\documentclass[pdflatex,12pt]{scrartcl}

\usepackage{array}
\usepackage{supertabular}
\usepackage{amsmath,amssymb}
\usepackage{bm}
\usepackage{graphicx}

\newcommand{\s}{\bm{\varSigma}}
\newcommand{\si}{\s^{-1}}
\newcommand{\lam}{\bm{\lambda}}
\newcommand{\Gam}{\bm{\varGamma}}
\newcommand{\diag}{\text{diag}}
\newcommand{\tr}{\text{tr}}
\newcommand{\xx}{\bm{x}}
\newcommand{\yy}{\bm{y}}
\newcommand{\el}{\bm{l}}
\newcommand{\uu}{\bm{u}}
\newcommand{\hh}{\bm{H}}
\newcommand{\wai}{\sum_{i}}

\newcommand{\waij}{\sum_{i<j}}
\newcommand{\tteta}{\Theta}

\newtheorem{prop}{Proposition}
\newtheorem{coro}{Corollary}

\DeclareMathOperator{\covae}{\stackrel{\mathit e}{\nabla}}
\DeclareMathOperator{\covam}{\stackrel{\mathit m}{\nabla}}
\DeclareMathOperator{\argmin}{argmin}

\makeatletter
  \def\mathcomposite{%
     \@ifstar
        {\def\@mathcomposite@option{%
            \baselineskip\z@skip\lineskiplimit-\maxdimen}%
         \@mathcomposite}%
        {\let\@mathcomposite@option\offinterlineskip
         \@mathcomposite}}
  \def\@mathcomposite{%
     \@ifnextchar[\@@mathcomposite{\@@mathcomposite[0]}}
  \def\@@mathcomposite[#1]#2#3#4{%
     #2{\mathchoice
        {\@mathcomposite@{#1}{#3}{#4}\displaystyle{1}}%
        {\@mathcomposite@{#1}{#3}{#4}\textstyle{1}}%
        {\@mathcomposite@{#1}{#3}{#4}%
         \scriptstyle\defaultscriptratio}%
        {\@mathcomposite@{#1}{#3}{#4}%
         \scriptscriptstyle\defaultscriptscriptratio}}}
  \def\@mathcomposite@#1#2#3#4#5{%
     \vcenter{\m@th\@mathcomposite@option
        \dimen@\f@size\p@\dimen@#1\dimen@\dimen@#5\dimen@
        \divide\dimen@ 18
        \edef\@mathcomposite@skipamount{\the\dimen@}%
        \ialign{\hfil$#4##$\hfil\cr
           #2\crcr
           \noalign{\vskip\@mathcomposite@skipamount}%
           #3\crcr}}}
  \makeatother

\title{Inference on the eigenvalues of the covariance matrix of a multivariate normal distribution--geometrical view--}
\author{Yo Sheena\thanks{Faculty of Economics, Shinshu University}}

\date{September 2012}

\begin{document}
\maketitle

\begin{abstract}
We consider inference on the eigenvalues of the covariance matrix of a multivariate normal distribution. The family of multivariate normal distributions with a fixed mean is seen as a Riemannian manifold with Fisher information metric. Two submanifolds naturally arises; one is the submanifold given by the fixed eigenvectors of the covariance matrix, the other is the one given by the fixed eigenvalues. We analyze the geometrical structures of these manifolds such as metric, embedding curvature under $e$-connection or $m$-connection. Based on these results, we study 1) the bias of the sample eigenvalues, 2) asymptotic variance of estimators, 3) the asymptotic information loss caused by neglecting the sample eigenvectors, 4) the derivation of a new estimator that is natural from a geometrical point of view.
\end{abstract}
\noindent
MSC(2010) \textit{  Subject Classification}: Primary 62H05; Secondary 62F12\\
\textit{t Key words and phrases:} curved exponential family, information loss, Fisher information metric, embedding curvature, affine connection, positive definite matrix.

\section{Introduction}
\label{section:int}
Consider a normal distribution with zero mean and an unknown covariance matrix, $N(\bm{0}, \s).$  Let denote the eigenvalues of $\s$ by 
$$
\lam=(\lambda_1, \ldots, \lambda_p),\quad \lambda_1> \ldots > \lambda_p
$$
and eigenvectors matrix by $\Gam$, hence we have the spectral decomposition
\begin{equation}
\label{spectdecomp_sigma}
\s=\Gam \bm{\Lambda} \Gam^t,  \quad \bm{\Lambda}=\diag(\lam),
\end{equation}
where $\diag(\lam)$ means the diagonal matrix with the $i$th diagonal element $\lambda_i$.
It is needless to say that the inference on $\s$ is an important task in many practical situations in such a diversity of fields as engineering, biology, chemistry, finance, psychology etc. Especially we often encounter the cases where the property of interest depends on $\s$ only through its eigenvalues $\lam$. We treat an inference problem on the eigenvalues $\lam$ from a geometrical point of view. 

Treating the  family of normal distributions $N(\bm{\mu},\s)$ ($\bm{\mu}$ is not necessarily zero) as a Riemmanian manifold has been done by several authors.  For example, see Fletcher and Joshi, \cite{Fletcher&Joshi}, Lenglet et al. \cite{Lenglet_et_al}, Skovgaard \cite{Skovgaard}, Smith \cite{Smith}, Yoshizawa and Tanabe \cite{Yoshizawa&Tanabe}. When $\mu$ euqals zero, the family of normal distributions $N(\bm{0},\s)$ can be taken as a manifold (say $\mathcal{S}$) with a single coordinate system $\s$. Hence, $\mathcal{S}$ is identified with the space of symmetric positive definite matrices. Geometrically analyzing the space of symmetric positive definite matrices has been an interesting topic in a mathematical or engineering point of view.  Refer to Moakher and Z\'{e}ra\"{i} \cite {Moakher&Zerai}, Ohara et al. \cite{Ohara_et_al} and Zhang et al. \cite{Zhang_et_al} as well as the above literature.

In this paper, we analyze $\mathcal{S}$ from the standpoint of information geometry while focusing on the inference on the eigenvalues of $\s.$ 
The paper is aimed to make a contribution in two regards: 1) The geometrical structure of $\mathcal{S}$ is analyzed in view of the eigenvalues and eigenvectors of $\s$; 2) Some statistical problems on the inference for $\lam$ are explained in the geometrical terms.

We summarize the inference problem for $\lam$.  Based on independent $n$ samples $\xx_i=(x_{i1},\ldots,x_{ip})',\ i=1,\ldots, n$ from  $N(\bm{0}, \s)$, we want to make inference on the unknown $\lam$. We confine ourselves to the classical case where $n\geq p$. It is well-known that the  product-sum matrix
$$
\bm{S}=\sum_{i=1}^n \xx_i \xx_i^t
$$
is sufficient statistic for both unknown $\lam$ and $\Gam$. 
The spectral decomposition of $\bm{S}$ is given by 
$$
\bm{S}=\hh\bm{L}\hh^t, \quad \bm{L}=\diag(\el),
$$
where 
$$
\el=(l_1,\ldots,l_p), \quad l_1>\ldots>l_p>0 \ a.e.
$$
are the eigenvalues of $\bm{S}$, and $\hh$ is the corresponding eigenvectors matrix. This decomposition gives us two statistics available, i.e. the sample eigenvalues $\el$ and the sample eigenvectors $\hh$.
However it is almost customary that we only use the sample eigenvalues, discarding the information contained in $\hh$. In the past literature on the inference for the population eigenvalues, every notable estimator is based simply on the sample eigenvalues. See Takemura \cite{Takemura}, Dey and Srinivasan \cite{Dey&Srinivasan},  Haff \cite{Haff}, Yang and Berger \cite{Yang&Berger} for orthogonally invariant estimators of $\s$; Dey \cite{Dey},  Hydorn and Muirhead   \cite{Hydorn&Muirhead}, Jin \cite{Jin}, Sheena and Takemura \cite{Sheena&Takemura} for direct estimators of $\lam$. Since we do not have enough space to state the concrete form of each estimator, we just mention Stein's estimator as a pioneering work for "shrinkage" estimator of $\s$. In general, an orthogonally invariant estimator of $\s$ is given by 
\begin{equation}
\label{orthogonally_inv_est}
\hat\s=\hh\bm{\Phi}\hh^t, \quad \bm{\Phi}=\diag(\phi_1(\el),\ldots,\phi_p(\el)).
\end{equation}
The estimator of $\lam$ is given by the eigenvalues of $\hat\s$, that is, $(\phi_1(\el),\ldots,\phi_p(\el))$.
The sample covariance matrix (M.L.E. estimator) $\bar{\bm{S}}\triangleq n^{-1}\bm{S}$ gives the estimator of $\lam$ as $\phi_i(\el)=n^{-1}l_i,\ i=1,\ldots,p$, while Stein's "shrinkage" estimator gives birth to 
\begin{equation}
\label{stein's_est}
\phi_i(\el)=l_i/(n+p+1-2i),\quad i=1,\ldots,p. 
\end{equation}
Stein's estimator assigns the lighter (heavier) weight to the larger (smaller) sample eigenvalues, hence the diversity of $\el$ is shrunk. This estimator is quite simple and performs much better than M.L.E. (see \cite{Dey&Srinivasan} ). Unlike Stein's estimator, many estimators in the above literature are not explicitly given or too complicated for immediate use. Nonetheless they all have one common feature. The derived estimators of $\lam$ only depends on $\el$.

In a sense it is natural to implicitly associate the sample eigenvalues  to the population eigenvalues, and the sample eigenvectors to the population counterpart. However the sample eigenvalues are not sufficient for the unknown population eigenvalues. Therefore it is important to evaluate how much information is lost by neglecting the sample eigenvectors.  Following  Amari \cite{Amari1}, we gain an understanding of the asymptotic information loss with geometric terms such as Fisher information metric and embedding curvatures.

Another statistically interesting topic is the bias of $n^{-1}\el$. It is well known that $n^{-1}\el$ is largely biased and the estimators mentioned above are all modification of $n^{-1}\el$ to correct the bias, that is, "shrinkage estimators." We show that the bias is closely related to the embedding curvatures. Moreover the geometric structure of $\mathcal{S}$ naturally leads us to a new estimator, which is also a shrinkage estimator.

The organization of this paper is as follows: In the former part (Section \ref{section:riemannian_manifold_and_metric} and Section \ref{section:embedding_curvature}), we describe the geometrical structure of $\mathcal{S}$ in view of the spectral decomposition \eqref {spectdecomp_sigma}. In Section \ref{section:riemannian_manifold_and_metric}, we observe $\mathcal{S}$ as a Riemannian manifold endowed with Fisher information metrics. In Section \ref{section:embedding_curvature}, we treat two submanifolds of $\mathcal{S}$, a submanifold given by the fixed eigenvectors and the one given by the fixed eigenvalues. The embedding curvatures of these submanifolds are explicitly given. We will show  that the bias of $\el$ is closely related to the curvatures. In the latter part (Section \ref{section:estimation of lambda when Gamma is known} and \ref{section:estimation of lambda when Gamma is unknown}), we consider the estimation problem of $\lam$. 
In Section \ref{section:estimation of lambda when Gamma is known}, we describe the asymptotic variance of estimators when $\Gam$ is known (Section \ref{subsec:asympto variance}) and the asymptotic information loss caused by discarding the sample eigenvectors $\hh$ (Section \ref{subsec:info loss}). The asymptotic information loss could be measured by the difference in the asymptotic variance between two certain estimators. In Section \ref{section:estimation of lambda when Gamma is unknown} for the case when $\Gam$ is unknown, we propose a new estimator of $\lam$, which is naturally derived from a geometric point of view. In the last section, some comments are made for further research.  All the proofs  are collected in Appendix.  

Unfortunately we do not have enough space to explain the geometrical concepts used in this paper. Please refer to  Boothby \cite{Boothby}, Amari \cite{Amari2}, Amari and Nagaoka \cite{Amari&Nagaoka}. 
\section{Riemannian Manifold and Metric}
\label{section:riemannian_manifold_and_metric}

The density of the normal distribution $N(\bm{0},\s)$ is given by
$$
f_{\s}(\xx)=(2\pi)^{-p/2} |\s|^{-1/2} \exp\left(-\frac{1}{2} \xx^t \si \xx \right),\quad \xx=(x_1,\ldots, x_p) \in R^p
$$
If we let $\sigma_{ij}$ and $\sigma^{ij}$ denote the $(i,j)$ element of respectively $\s$ and $\si$, then the log likelihood equals
\begin{equation}
\label{lognormal}
\begin{split}
\log f_{\s}(\xx)&=\wai x_i^2 \left(-\sigma^{ii}/2\right)+\waij x_i x_j \left(-\sigma^{ij} \right)-(p/2)\log 2\pi-(1/2)\log|\s| \\
            &=\wai y_{ii}\theta^{ii} + \waij y_{ij} \theta^{ij} -\psi(\tteta) \ (\text{say }l(\yy;\tteta)),
\end{split}
\end{equation}
where $\tteta=(\theta^{ij})_{i \leq j}$ and $\yy=(y_{ij})_{i\leq j}$ are given by
\begin{equation}
\label{natural_parameters}
\left\{
\begin{aligned}
\theta^{ii}&=(-1/2)\sigma^{ii},\qquad &&i=1,\ldots,p, \\
\theta^{ij}&=-\sigma^{ij},\qquad &&1\leq i <j  \leq p, \\
y_{ii}&=x_{i}^2,\qquad &&i=1,\ldots,p, \\
y_{ij}&=x_i x_j,\qquad &&1\leq i <j  \leq p, 
\end{aligned}
\right.
\end{equation}
and 
\begin{equation}
\psi(\tteta)=(p/2)\log 2\pi +(1/2) \log |\s(\tteta)|.
\end{equation}
The summations $\Sigma_{i}$, $\Sigma_{i<j}$ in the equation \eqref{lognormal} are abbreviations respectively for $\sum_{i=1}^p$ and $\sum_{1\leq i <j \leq p}$, and we will use these kinds of notations implicitly hereafter.

The expression \eqref{lognormal} gives natural coordinate system $\tteta$ of the manifold $\mathcal{S}$ as a full exponential family.  Another coordinate system, so called expectation parameters, is also useful, which is defined as;
\begin{equation}
\label{expectation_parameters}
\sigma_{ij}=E(y_{ij}),\qquad 1\leq i \leq j \leq p.
\end{equation}

For the analysis of the information carried by $\el$ and $\hh$, we need to prepare another coordinate system. The matrix exponential expression of an orthogonal matrix $\bm{O}$ is given by
\begin{equation}
\bm{O}=\exp\bm{U}=\bm{I}_p+\bm{U}+\frac{1}{2}\bm{U}^2+\frac{1}{3!}\bm{U}^3+\cdots,
\end{equation}
where $\bm{I}_p$ is the $p$-dimensional unit matrix,  $\bm{U}$ is a skew-symmetric matrix and parametrized by $\uu=(u_{ij})_{1\leq i<j \leq p}$ as
$$
(\bm{U})_{ij}= \left\{
\begin{aligned}
u_{ij},     &&\text{if }  1\leq i <j \leq p,\\
-u_{ij},    &&\text{if }  1\leq j <i \leq p,\\
0,         &&\text{if } 1\leq i=j \leq p.
\end{aligned}
\right.
$$
The function $\exp\bm{U}$ is diffeomorphic, and $\uu$ gives  "normal coordinate" for the group of orthogonal matrices (see  (6.7) in Boothby\cite{Boothby} or Th. A9.11 of Muirhead\cite{Muirhead}). We can use this coordinate as local system around $\bm{I}_p$ and construct an atlas for the entire space of p-dimensional orthogonal matrices (note this space is compact); for each $\Gam$, there exists an open neighborhood  and some open ball $B$ in $R^{p(p-1)/2}$ around the origin such that these spaces are diffeomorphic by the function $\Gam \exp \bm{U}(\uu)$ on $B$.

We will use $(\lam, \uu)$ as the third coordinate system of $\mathcal{S}$ and call it "spectral coordinate (system)". Notice that this coordinate system is associated with the following submanifolds in $\mathcal{S}$. If we fix $\Gam$ in \eqref{spectdecomp_sigma}, then we get a submanifold $\mathcal{M}(\Gam)$ embedded in $\mathcal{S}$ with a coordinate system $\lam$.  This is a subfamily in $N(\bm{0},\s)$ and called curved exponential family. Its log-likelihood is expressed, as we emphasize it as a function of $\lam$, to be
\begin{equation}
l(\yy;\tteta(\lam))=\wai y_{ii} \theta^{ii}(\lam)+\waij y_{ij}\theta^{ij}(\lam)-\psi(\tteta(\lam)).
\end{equation}
On the contrary, if we fix $\lam$ in \eqref{spectdecomp_sigma}, we get another submanifold $\mathcal{A}(\lam)$ in $\mathcal{S}$, whose coordinate system is given by $\uu$ in a neighborhood of each point of $\mathcal{A}(\lam)$.  Its log-likelihood expression is given by
\begin{equation}
l(\yy;\tteta(\uu))=\wai y_{ii} \theta^{ii}(\uu)+\waij y_{ij}\theta^{ij}(\uu)-\psi(\tteta(\uu)).
\end{equation}

First we consider a metric, that is, a field of symmetric, positive definite, bilinear form on $\mathcal{S}$. The statistically most natural metric is Fisher information metric. Suppose $\{f(\xx;\bm{\theta})\}$ is a  parametric family of probability density functions, whose coordinate as a manifold is given by $\bm{\theta}=(\theta_1,\ldots,\theta_p)$. Then the $(i,j)$ component of Fisher information metric with respect to $\bm{\theta}$ is  given by 
$$
E_{\bm{\theta}}\left[\frac{\partial}{\partial \theta_i}\log f(\xx; \theta) \frac{\partial}{\partial \theta_j}\log f(\xx; \theta)\right].
$$
For the multivariate normal distribution family, $N(\bm{\mu},\s)$ ($\bm{\mu}$, the mean parameter is also included),  Skovgaard \cite{Skovgaard} gives a clear form of Fisher information metric. The  tangent vector space at a fixed point $\s$ w.r.t. $(\sigma_{ij})_{i\le j}$ coordinate can be identified with the space of symmetric matrices. For any symmetric matrix $\bm{A}$, $\bm{B}$, the metric with respect to the $\s=(\sigma_{ij})$ coordinate system is given by
\begin{equation}
\label{metric_simple_form}
\frac{1}{2} \text{tr}\Bigl(\si \bm{A} \si \bm{B}\Bigr).
\end{equation}

We are interested in Fisher information metric with respect to the spectral coordinate $(\lam,\uu)$. 
Let $\partial_a,\ \partial_b,\ \cdots$ denote the tangent vectors w.r.t. the $\lam$ coordinate, $\partial_{(s,t)},\ \partial_{(u,v)},\cdots$ denote the tangent vectors w.r.t. the $\uu$ coordinate. Namely 
$$
\partial_a\triangleq \frac{\partial }{\partial \lambda_a},\qquad  \partial_{(s,t)}\triangleq\frac{\partial }{\partial u_{st}}.
$$
These tangent vectors (exactly speaking, vector fields) are invariant with respect to any orthogonal transformation of $\s$; For some orthogonal matrix $\bm{O}$, an orthogonal transformation $F$of $\mathcal{S}$ is defined as
\begin{equation}
\label{orthogotranse}
F(\s)=\bm{O}\s\bm{O}^t
\end{equation}
For any $\bm{O}, $
\begin{align}
F_*(\partial_a)&=\partial_a,\qquad &&1\leq a \leq p, \\
F_*(\partial_{(s,t)})&=\partial_{(s,t)},\qquad && 1\leq s<t \leq p,
\end{align}
where $F_*$ is the derivative of $F$.
\begin{prop}
\label{metric_spectra}
Let  $\langle\ , \ \rangle$ denote Fisher information metric based on $\xx \sim N(\bm{0},\s)$, then
the components of the metric with respect to $(\lam, \uu)$ is given as follows;
\begin{align*}
&g_{ab}\triangleq\langle\partial_a,\partial_b \rangle=(1/2)\lambda_a^{-2} \;\delta\bigl(a=b\bigr) \quad && 1\leq a, b \leq p, \\
&g_{a(s,t)}\triangleq \langle \partial_a,\partial_{(s,t)}\rangle =0 \quad && 1\leq a \leq p, \ 1\leq s < t \leq p,\\
&g_{(s,t)(u,v)}\triangleq \langle \partial_{(s,t)},\partial_{(u,v)}\rangle \\
&=(\lambda_s-\lambda_t)^2\lambda_s^{-1}\lambda_t^{-1}\;\delta\bigl((s,t)=(u,v)\bigr) &&1\leq s <t \leq p,\  1 \leq u < v \leq p.
\end{align*}
$\delta(\cdot)$ equals one if the logic inside the parenthesis is correct, otherwise zero.
\end{prop}

There are two remarkable properties of the metric for the spectral coordinate. First note that since the metric components matrix is diagonal,  $(\lam,\uu)$ is an orthogonal coordinate system, especially that the submanifolds $\mathcal{M}(\Gam)$ and $\mathcal{A}(\lam)$ are orthogonal to each other for any $\lam$ and $\Gam$. Second it is independent of $\Gam$, hence the metric stays constant with respect to the orthogonal transformation $F$ in \eqref{orthogotranse} for any orthogonal matrix $\bm{O}.$ (Second property is instantly derived from the expression \eqref{metric_simple_form}.) 

Theoretically, other metrics could be naturally implemented. Calvo and Oller \cite{Calvo&Oller} introduced Sigel metric. Lovri\'c et al. \cite{Lovric_et_al} considered the natural invariant metric from the standpoint of Riemannian symmetric space. The concrete forms of the both metrics are given by (3.4) and (3.2) in \cite{Lovric_et_al}. (The information metric \eqref{metric_simple_form} corresponds to (3.3) in \cite{Lovric_et_al}. See also Theorem 1 of Zhang \cite{Zhang_et_al}. )

Once a metric is given on the manifold $\mathcal{S}$,  a connection is needed for further geometrical analysis. Connection is an important "rule" which defines how a tangent space is shifted with an infinitesimal move in a differential manifold. Although connection has an infinite variation, the most commonly used one is  Levi-Civita connection. It is characterized as a unique torsion-free, metric-preserving connection. This connection is essential to consider a distance function on the manifold. Skovgaard \cite{Skovgaard} , Calvo and Oller \cite{Calvo&Oller}, Fletcher and Joshi \cite{Fletcher&Joshi}, Lenglet et al. \cite{Lenglet_et_al}, Lovri\'c et al. \cite{Lovric_et_al}, Moakhaer and Z\'{e}ra\"{i} \cite{Moakher&Zerai} analyze the manifold of the normal distributions under Levi-Civita connection. 

On the other hand, Amari [1] showed that "$\alpha$-connection" is suitable for statistical manifolds in general. He also found that e-connection ($\alpha=1$) and m-connection ($\alpha=-1$) are especially important for the asymptotic analysis of information loss for a curved exponential family. Amari and Kumon \cite{Amari&Kumon}, Kumon, Amari \cite{Kumon&Amari} and Eguchi \cite{Eguchi} gave further development along this line. Specifically in the relation with the multivariate normal distribution or $\mathcal{S}$, Ohara et al. \cite{Ohara_et_al}, Yoshizawa and Tanabe \cite{Yoshizawa&Tanabe} and Zhang et al. \cite{Zhang_et_al} considered the dual geometry ($\alpha$ and $-\alpha$ connections) of the manifolds.  Notice that Levi-Civita connection is $0$-connection and the "mean" between $e$-connection and $m$-connection. Therefore, using the results on geometric properties of $\mathcal{S}$ under $e$-connection and $m$-connection, we could also derive those under Levi-Civita connection. 

Since this paper is aimed for the statistical  inference on $\s$, we adopt $\alpha$-connections, especially $e$- and $m$-connections, hereafter. 
We conclude this section by mentioning the important fact that $\mathcal{S}$ is $e$-flat and $m$-flat, and corresponding affine coordinates are given respectively by $(\sigma^{ij})$ and $(\sigma_{ij})$.
\section{Embedding Curvatures}
\label{section:embedding_curvature}

Curvature, which is important property for an geometrical analysis, is defined based on a given  connection. A submanifold has both intrinsic and extrinsic curvatures. The latter describes how the submanifold is placed in the whole manifold, and called an embedding curvature or the second fundamental form. (The first fundamental form is the metric.)

In this section, we observe the embedding curvatures of $\mathcal{M}$ and $\mathcal{A}$ for the analysis of the distribution $(\el, \hh)$. Specifically we consider the following embedding curvatures;

1. Embedding curvature of $\mathcal{M}$ with respect to $e$-connection or $m$-connection. Its components w.r.t the spectral coordinate are given by
\begin{equation}
\stackrel{e}{H}_{ab(s,t)}\triangleq \langle \covae_{\partial_a}\!\!\partial_b\,, \partial_{(s,t)}\rangle, 
\qquad \stackrel{m}{H}_{ab(s,t)}\triangleq \langle \covam_{\partial_a}\!\!\partial_b\,, \partial_{(s,t)}\rangle, 
\end{equation}
where $\covae_{\partial_a}\!\!\partial_b$ is the covariant derivative of $\partial_b$ in the direction of $\partial_a$ with respect to $e$-connection. $\covam_{\partial_a}\!\!\partial_b$ is similarly defined.

2. Embedding curvature of $\mathcal{A}$ with respect to $m$-connection. Its components w.r.t the spectral coordinate are given by
\begin{equation}
\stackrel{m}{H}_{(s,t)(u,v)a}\triangleq\langle\covam_{\partial_{(s,t)}}\!\!\partial_{(u,v)}\,,\partial_{a}\rangle,
\end{equation}
where $\covam_{\partial_{(s,t)}}\!\!\partial_{(u,v)}$ is the covariant derivative of $\partial_{(s,t)}$ in the direction of $\partial_{(u,v)}$ with respect to $m$-connection.

On these curvatures at the point $(\lam, \Gam)$, we have the following results.
\begin{prop}
\label{embedding_c}
For  $1\leq a,b \leq p,\ 1\leq s<t \leq p, $
\begin{equation}
\label{embedding_c_H}
\stackrel{e}{H}_{ab(s,t)}=\stackrel{m}{H}_{ab(s,t)}=0.
\end{equation}
For $1\leq a \leq p, \ 1\leq s<t \leq p,\ 1\leq u<v \leq p, $
\begin{equation}
\label{embedding_curvature_A}
\stackrel{m}{H}_{(s,t)(u,v)a}=
\begin{cases}
\lambda_a^{-2}(\lambda_t-\lambda_a),  \text{ if $s=u=a$, $t=v$,}\\
\lambda_a^{-2}(\lambda_s-\lambda_a),  \text{ if $s=u$, $t=v=a$,}\\
0, \text{ otherwise.}
\end{cases}
\end{equation}
\end{prop}

Another expression of the embedding curvature of $\mathcal{A}$ is given by
\begin{equation}
\label{embedd_m_curve_A}
\stackrel{m}{H^{a}}_{\!\!\!\!(s,t)(u,v)}\triangleq \sum_{b}\stackrel{m}{H}_{(s,t)(u,v)b}g^{ba},
\end{equation}
With this notation, the orthogonal projection of the covariant derivative
$$
\covam_{\partial_{(s,t)}}\!\!\partial_{(u,v)}
$$
onto the tangent space of $\mathcal{M}$ is given by
$$
\sum_{a}\stackrel{m}{H^{a}}_{\!\!\!\!(s,t)(u,v)}\partial_a.
$$
From Proposition \ref{metric_spectra}, \ref{embedding_c}, we have
\begin{equation}
\label{embedding_c_2}
\stackrel{m}{H^{a}}_{\!\!\!\!(s,t)(u,v)}=2(\lambda_t-\lambda_a)\delta(s=u=a,\,t=v)+2(\lambda_s-\lambda_a)\delta(s=u,\,t=v=a),
\end{equation}
hence
$$
\sum_{a}\stackrel{m}{H^{a}}_{\!\!\!\!(s,t)(u,v)}\partial_a=
\begin{cases}
2(\lambda_t-\lambda_s)\partial_s +2(\lambda_s-\lambda_t)\partial_t, &\text{if $(s,t)=(u,v),$}\\
0, &\text{otherwise.}
\end{cases}
$$
Similarly another embedding curvature components $\mathcomposite[-2]{\mathrel}{\hspace{-5mm}{\scriptsize\textit e}}{H^{(s,t)}_{ce}}$ is defined as 

\begin{equation}
\label{compo_e-embedd_C_M}
\mathcomposite[-2]{\mathrel}{\hspace{-5mm}{\scriptsize\textit e}
}{H^{(s,t)}_{ab}}=\sum_{u<v}\stackrel{e}{H}_{ab(u,v)}g^{(u,v)(s,t)}
\end{equation}
and actually it vanishes
\begin{equation}
\label{e-embedding_c_M_upper_comp_vanish}
\mathcomposite[-2]{\mathrel}{\hspace{-5mm}{\scriptsize\textit e}
}{H^{(s,t)}_{ab}}=0,\quad
1\leq a,b \leq p,\ 1\leq s<t \leq p.
\end{equation}

An embedding curvature has full information about the "extrinsic curvature" of the embedded submanifold in any direction.  Sometimes it is convenient to compress it into a scalar measure of the curvature. "Statistical curvature" by Efron (see Efron \cite{Efron}, Murray and Rice \cite{Murray&Rice}) is such a measure; For $\mathcal{A}$, it is defined by (see p.159 of Amari \cite{Amari2}) 
$$
\gamma({\mathcal{A}})\triangleq\sum_{1\leq a, b \leq p}\;\sum_{s<t,u<v,o<p,q<r}\stackrel{m}{H}_{(s,t)(u,v)a}\;\stackrel{m}{H}_{(o,p)(q,r)b}\;g^{(s,t)(o,p)}\;g^{(u,v)(q,r)}\;g^{ab},
$$
which attains the following value at the point $(\lam, \Gam)$.
\begin{coro}
\label{efron_C}
$$
\gamma({\mathcal{A}})=2\sum_{a<b} \frac{\lambda_a^2+\lambda_b^2}{(\lambda_a-\lambda_b)^2}
$$
\end{coro}

From these results, we notice that if $\mathcal{S}$ is endowed with $m$-connection, then 1) the embedding curvatures and the statistical curvatures of $\mathcal{A}$ are independent of $\Gam$,  2) any one-parameter curve $(\lam, \Gam(\bm{u}))$ given by a parameter $u_{(s,t)},\ s<t$, where $\lam$ and the other elements of $\bm{u}$ are fixed, is curved in the direction of $\partial_t-\partial_s$ and contained in a two-dimensional plane composed by $\partial_{(s,t)}$ and $\partial_t-\partial_s$, 3) the statistical curvature of $\mathcal{A}$ could be quite large when $\lam$ are close to each other, while $\mathcal{M}$ is flat everywhere. 

Here we introduce another submanifold $\tilde{\mathcal{A}}$ which is contrasting to $\mathcal{A}$ in the sense that $\tilde{\mathcal{A}}$ is flat with respect to $m$-connection. 
For a point $(\lam, \Gam)$,  let
$$
\tilde{\mathcal{A}}(\lam, \Gam)\triangleq\{\s \in \mathcal{S} \;|\; (\Gam^t\s\Gam)_{ii}=\lambda_{i},\ 1\leq \forall i \leq p\}.
$$
We easily notice that $\tilde{\mathcal{A}}$ is the minimum distance points with respect to Kullback-Leibler divergence.  That is,
$$
\tilde{\mathcal{A}}(\lam, \Gam)=\{ \s \in \mathcal{S} \;|\; \argmin_{\tilde{\lam}} KL(\s, \Gam {\textrm diag}(\tilde{\lambda}_1,\ldots,\tilde{\lambda}_p)\Gam^t)=\lam\},
$$
where $KL(\s, \tilde\s)$ is the Kullback-Leibler divergence between $N(\bm{0},\s)$ and $N(\bm{0},\tilde\s)$, which is specifically given by
$$
{\textrm tr}(\s \tilde{\s}^{-1})-\log|\s \tilde{\s}^{-1}|-p.
$$
The minimum distance points with respect to the Kullback-Leibler divergence consists of all the points on the $m$-geodesics which pass through the point $(\lam, \Gam)$ and are orthogonal to $\mathcal{M}(\Gam)$ at that point. (See Theorem in A2 of Amari [1]).

\begin{figure}[t]
\centering
\includegraphics[width=8cm,bb=9 9 458 558]{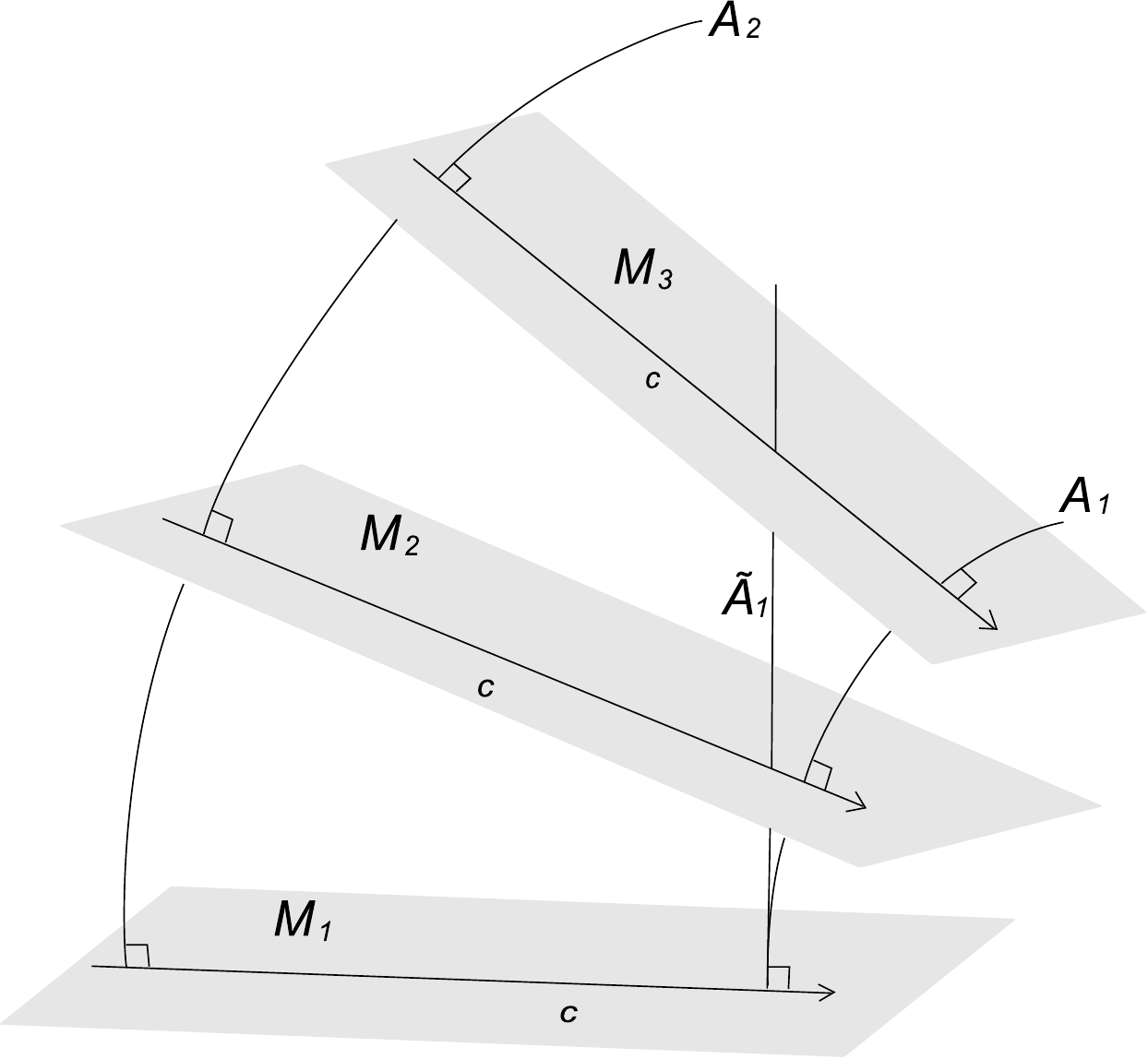}
\caption{Submanifolds of $\mathcal{S}$ when $p=2$, $\mathcal{M}$, $\mathcal{A}$ and $\tilde{\mathcal{A}}$}
\label{fig1}
\end{figure}

We can visualize the structure of $\mathcal{S}$ endowed with $m$-connection for the two dimensional case. See Figure \ref{fig1}, where $\mathcal{M}_i\triangleq\mathcal{M}(\Gam_i),\ i=1,\ldots, 3$, $\mathcal{A}_i\triangleq \mathcal{A}(\lam_i),\ i=1,2$ and $\tilde{\mathcal{A}}_1\triangleq \tilde{\mathcal{A}}(\lam_1,\Gam_1)$ are drawn.  When $p=2$, $\mathcal{M}$ is a two-dimensional autoparallel submanifold with the affine coordinate $(\lambda_1, \lambda_2)$, while $\mathcal{A}$ is a one-dimensional submanifold with an coordinate $u_{(1,2)}$. As it is seen in Proposition \ref{metric_spectra}, all the tangent vectors $\partial_1(\triangleq\frac{\partial}{\partial\lambda_1})$, $\partial_2(\triangleq\frac{\partial}{\partial\lambda_2})$, $\partial_{(1,2)}(\triangleq\frac{\partial}{\partial u_{(1,2)}})$ are orthogonal to each other. $\tilde{\mathcal{A}}$ is a "straight" line which is also orthogonal to $\mathcal{M}$.
The arrow on $\mathcal{M}$ is the line $\{\lam | \lambda_1+\lambda_2\text{ is constant}\}$, and the arrow head indicates the direction in which $c\triangleq\lambda_2/\lambda_1$ increases. The statistical curvature turns out to be the increasing function of $c$ ;
$$
\gamma(\mathcal{A})=2 \frac{1+c^2}{(1-c)^2}.
$$

\begin{figure}[t]
\centering
\includegraphics[width=10cm,bb=9 9 458 243]{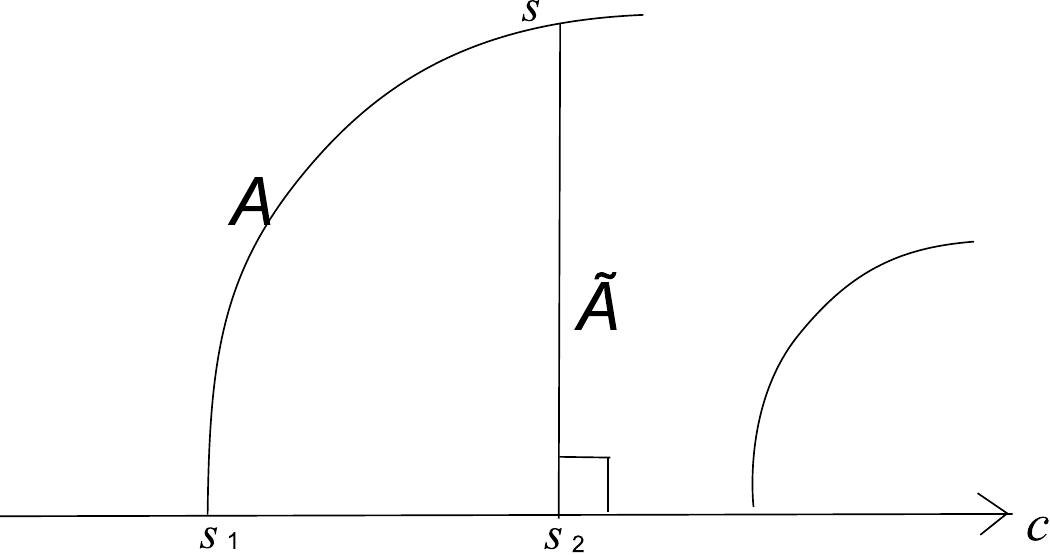}
\caption{ Horizontal perspective of $\mathcal{A}$ and $\tilde{\mathcal{A}}$ on the plane $\mathcal{M}$ when $p=2$ }
\label{fig2}
\end{figure}

We can analyze the bias of $\bar{l}_i\triangleq n^{-1}l_i,\ i=1,\ldots, p$ from the  geometrical structure of $\mathcal{S}$. It is well known that $E[\,\bar{l}_i\,]\ (i=1,\ldots,p)$ majorizes $\lambda_i\ (i=1,\ldots,p)$, that is, 
\begin{equation}
\label{bias_l}
\sum_{i=1}^j E[\,\bar{l}_i\,] \geq \sum_{i=1}^j \lambda_i,\quad 1\leq \forall j \leq p-1, \qquad \sum_{i=1}^p E[\,\bar{l}_i\,] = \sum_{i=1}^p \lambda_i.
\end{equation}
The bias $E[\,\bar{l}_i\,]$ is quite large when $n$ is small and $\lambda_i$'s are close to each other (see Lawley \cite{Lawley}, Anderson \cite{Anderson}). For the case $p=2$, 
\begin{equation}
\label{bias_l_p=2}
E[\,\bar{l}_1\,] \geq \lambda_1,\qquad E[\,\bar{l}_2\,] \leq \lambda_2,\qquad E[\,\bar{l}_1\,]+E[\,\bar{l}_2\,]=\lambda_1+\lambda_2.
\end{equation}
Suppose a sample $\bar{\bm{S}}\triangleq n^{-1}\bm{S}$ takes the value at a point $s \in \mathcal{S}.$ Let $s_1$ denote the point  on $\mathcal{M}(\Gam)$ designated by the  eigenvalues of $\bar{\bm{S}}$, namely  $\bar{\bm{l}}\triangleq(\bar{l}_1,\bar{l}_2)$. The curve $\mathcal{A}(\bar{\bm{l}})$ connects $s$ and $s_1$. If we define $s_2$ as the point on $\mathcal{M}(\Gam)$ designated by $\hat{\lam}\triangleq(\hat{\lambda}_1,\hat{\lambda}_2)\triangleq((\Gam^t\bar{\bm{S}}\Gam)_{11},(\Gam^t\bar{\bm{S}}\Gam)_{22})$, then $\tilde{\mathcal{A}}(\hat{\lam},\Gam)$ connects $s$ and $s_2$. The three points $s$, $s_1$ and $s_2$ are on the same plane, and if we move from $s_1$ in the direction to $s_2$, then the statistical curvature of $\mathcal{A}$ increases (see Figure \ref{fig2}). 
 If we estimate $(\lambda_1,\lambda_2)$ by $\bar{\el}$, then the estimate is the point $s_1$, while for the unbiased estimator $\hat{\lam}$, the estimate is the point $s_2$. Since the $c$-coordinate of $s_1$ is always smaller than that of $s_2$, the estimator $(\bar{l}_1,\bar{l}_2)$ is likely to estimate $\lambda_1$ and $\lambda_2$ too apart, which causes the bias \eqref{bias_l_p=2}. It is also seen that the bias gets larger when $c$ approaches to one, that is, $\lambda_1$ and $\lambda_2$ get closer to each other.

Though the exact magnitude of the bias $E(\bar{l}_a)-\lambda_a$ is hard to evaluate, the asymptotic bias can be evaluated. This can be also described with embedding curvatures (see (5.4) of Amari \cite{Amari2});
$$
E(\bar{l}_a-\lambda_a)=-\frac{1}{2n}C^a+O(n^{-3/2}),
$$
where
$$
C^a=\sum_{c,d}\stackrel{m}{\varGamma^{a}}_{\!\!\!\!cd}g^{cd}+\sum_{s<t, u<v}\stackrel{m}{H^{a}}_{\!\!\!\!(s,t)(u,v)}g^{(s,t)(u,v)},\\
$$
and $\stackrel{m}{\varGamma^{a}}_{\!\!\!\!cd}$ is a $m$-connection coefficients of $\mathcal{M}$, which is defined by
\begin{equation}
\label{m-connection-coeff-M}
\stackrel{m}{\varGamma^{a}}_{\!\!\!\!cd\;}=\stackrel{m}{\varGamma}_{\!\!cdb}g^{ba},\quad
\stackrel{m}{\varGamma}_{\!\!cdb}\triangleq \langle \covam_{\partial_c}\!\!\partial_d\,, \partial_{b}\rangle.
\end{equation}
Since $\mathcal{M}$ is autoparallel in $m$-flat $\mathcal{S}$, 
\begin{equation}
\label{conection_coeff_M_vanish}
\stackrel{m}{\varGamma^{a}}_{\!\!\!\!cd}=\stackrel{m}{\varGamma}_{\!\!cdb}=0, \quad 1\leq a,b,c,d \leq p.
\end{equation} Hence we have the following equation from Proposition \ref{metric_spectra} and \eqref{embedding_c_2}.
\begin{align}
\label{asympto_bias}
C^a(\lam)&=\sum_{a<t}\stackrel{m}{H^{a}}_{\!\!\!\!(a,t)(a,t)}g^{(a,t)(a,t)}+\sum_{s<a}\stackrel{m}{H^{a}}_{\!\!\!\!(s,a)(s,a)}g^{(s,a)(s,a)}\nonumber\\
     &=2\sum_{t\ne a} \frac{\lambda_a \lambda_t}{\lambda_t-\lambda_a}.
\end{align}
This bias was originally derived by the perturbation method in  Lawley \cite{Lawley}.
\section{Estimation of $\lam$ when $\Gam$ is known}
\label{section:estimation of lambda when Gamma is known}

We consider an estimation problem when  $\Gam$ is known to be $\Gam^0$. 
From a practical point of view, the case when $\Gam$ is known is not of much interest compared to the general case where both $\Gam$ and $\lam$ are unknown. However as we will show in this section, the asymptotic information loss caused by discarding the sample eigenvectors (Section \ref{subsec:info loss}) are closely related to the asymptotic variance difference between two certain estimators  (Section \ref{subsec:asympto variance}).  Both asymptotic variance and information loss are described with geometrical terms.
\subsection{Asymptotic variance of the estimators of $\lam$}
\label{subsec:asympto variance}
In a general term, the subfamily (submanifold) $\mathcal{M}(\Gam^0) (\triangleq \{\s \in \mathcal{S} | \Gam(\s)= \Gam^0$\}) in $\mathcal{S}$ is a  "curved" exponential family, since it is a subfamily in an exponential family $\mathcal{S}$. In a usual case, a subfamily is not "flat", hence the term "curved" is used. However  as you can see from \eqref{embedding_c_H},  $\mathcal{M}(\Gam^0)$ is autoparallel in $m(e)$-flat $\mathcal{S}$, and intrinsically $m(e)$-flat (see e.g. Theorem 1.1 in \cite{Amari&Nagaoka}).

We are supposed to estimate unknown coordinate $\bm{\lambda}$ of  $\mathcal{M}(\Gam^0)$ using an estimator $\hat{\lam}=(\hat{\lambda}_1,\ldots,\hat{\lambda}_p)$ of some kind. An estimator $\hat{\lam}(\bm{S})$  is specified by its inverse image $\hat{\lam}^{-1}(\lam)$
\begin{equation}
\label{Inv_image_estimator}
\hat{\mathcal{A}}(\lam)\triangleq \hat{\lam}^{-1}(\lam)=\{\s \in \mathcal{S} \, | \, \hat{\lam}(\s)=\lam\}.
\end{equation}
This is another submanifold in $\mathcal{S}$, where we will use $\bm{u}$ as a coordinate system. 

A consistent estimator $\hat{\lam}$ is called first-order (Fisher) efficient if the first order term (i.e. $O(n^{-1}$) order term) w.r.t. the asymptotic expansion of the variance (covariance) in $n$ is minimized among all (regular) estimators. Correct the bias of the first-order efficient estimator $\hat{\lam}$ up to the term of order $n^{-1}$, and let it be denoted by $\hat{\lam}^*\triangleq(\hat{\lambda}^*_1,\ldots,\hat{\lambda}^*_p)$. Amari showed (see e.g. Theorem 4.4 in \cite{Amari&Nagaoka})  that its asymptotic variance can be  described by the geometrical properties such as the metric and the embedding curvatures of $\mathcal{M}(\Gam^0)$ and $\hat{\mathcal{A}}$ ; For $1\leq a,b \leq p$,
\begin{equation}
\label{asympto_variance}
E[(\hat{\lambda}^*_a-\lambda_a)(\hat{\lambda}^*_b-\lambda_b)]=\frac{1}{n}g^{ab}+\frac{1}{2n^2}\{(\varGamma^m_M)^{2ab}+2(H_M^e)^{2ab}+(H_{\hat A}^m)^{2ab}\}+O(n^{-3})
\end{equation}
where
\begin{align*}
(\varGamma^m_M)^{2ab}&=\sum_{c,d,e,f}\mathcomposite[3]{\mathrel}{\hspace{-1mm}{\scriptsize \textit m}
}{\varGamma^{a}_{cd}}\; \mathcomposite{\mathrel}{\hspace{-2mm}{\scriptsize \textit m}}{\varGamma^{b}_{ef}}  g^{ce}g^{df},\\
(H_M^e)^{2ab}&=\sum_{c,d,e,f, s<t, u<v}\mathcomposite{\mathrel}{\hspace{-5mm}{\scriptsize \textit e}}{H^{(s,t)}_{ce}} \:\mathcomposite[-3]{\mathrel}{\hspace{-6mm}{\scriptsize \textit e}}{H^{(u,v)}_{df}}
g_{(s,t)(u,v)}\:g^{cd}g^{ea}g^{fb},
\\
(H_{\hat A}^m)^{2ab}&=\sum_{s<t, u<v, o<p, q<r}\mathcomposite[3]{\mathrel}{\hspace{-11mm}{\scriptsize \textit m}
}{H^{a}_{(s,t)(u,v)}} \:\mathcomposite[1]{\mathrel}{\hspace{-11mm}{\scriptsize \textit m}}{H^{b}_{(o,p)(q,r)}} g^{(s,t)(o,p)}g^{(u,v)(q,r)},
\end{align*}
$\mathcomposite[3]{\mathrel}{\hspace{-1mm}{\scriptsize \textit m}}{\varGamma^{a}_{cd}}$ and 
$\mathcomposite[-3]{\mathrel}{\hspace{-5mm}{\scriptsize \textit e}}{H^{(s,t)}_{ce}}$ are already defined in the previous section as the connection coefficients (see \eqref{m-connection-coeff-M}) or the embedding curvature components (see \eqref{compo_e-embedd_C_M}) of $\mathcal{M}$. They are defined independently of the particular estimator. $\mathcomposite[3]{\mathrel}{\hspace{-11mm}{\scriptsize \textit m}}{H^{a}_{(s,t)(u,v)}}$ are the components of the embedding $m$-curvature of $\hat{\mathcal{A}}$, which differ among the estimators. 

We apply this formula to the following two estimators, $\bm{l}^*=(l_1^*, \ldots, l_p^*)$ and  $\hat{\lam}=(\hat{\lambda}_1,\ldots,\hat{\lambda}_p).$ The former is the bias-corrected sample eigenvalues, which is given, using \eqref{asympto_bias}, by 
\begin{equation}
\label{bias_correct_esti}
l_a^*=\bar{l}_a+\frac{1}{2n}C^a(\el)=\bar{l}_a+\frac{1}{n}\sum_{t\ne a}\frac{l_a l_t}{l_t-l_a},\quad a=1,\ldots,p,
\end{equation}
and the latter is defined by
\begin{equation}
\label{m.l.e.}
\hat{\lambda}_a=((\Gam^0)^t\bar{\bm{S}}\Gam^0)_{aa},\quad a=1,\ldots,p,
\end{equation}
which is  (exactly) unbiased. In fact $\hat{\lam}$ is the maximum likelihood estimator for the case $\Gam$ is known. Notice that for  $\bm{l}$, $\hat{\mathcal{A}}(\lam)=\mathcal{A}(\lam)$ and that for $\hat{\lam}$, $\hat{\mathcal{A}}(\lam)=\tilde{\mathcal{A}}(\lam, \Gam^0)$. The first-order efficiency of the both estimators are guaranteed by the orthogonality to $\mathcal{M}(\Gam^0)$ of $\mathcal{A}(\lam)$ and $\tilde{\mathcal{A}}(\lam, \Gam^0)$.

The terms $(\varGamma^m_M)^{2ab}$ and $(H_M^e)^{2ab}$, which are related to the submanifold $\mathcal{M}$, hence common to the both estimators, vanish, because of 
\eqref{e-embedding_c_M_upper_comp_vanish} and \eqref{conection_coeff_M_vanish}.
The term $(H_{\hat A}^m)^{2ab}$ is different between the two estimators.  As we observed in the previous section, $\mathcal{A}(\lam)$ is not autoparallel in $\mathcal{S}$ (see \eqref{embedding_curvature_A} ). On the other hand, $\tilde{\mathcal{A}}(\lam, \Gam^0)$ is autoparallel  in $\mathcal{S}$, hence $(H_{\hat A}^m)^{2ab}$ vanishes. Consequently the following results are gained.
\begin{prop}
\label{asymp_var_result}
For $1\leq a, b \leq p$,
\begin{align}
&E[(l_a^*-\lambda_a)(l_b^*-\lambda_b)]  (\triangleq V^{ab}(\el^*)) \nonumber\\
\label{asymp_var_result_sample_eigen}
&=
\begin{cases}
\displaystyle
\frac{2}{n}\lambda_a^2+\frac{2}{n^2}\sum_{t\ne a}\frac{\lambda_a^2 \lambda_t^2}{(\lambda_t-\lambda_a)^2}+O(n^{-3}), & \text{if $a=b$,}\\
&\\
\displaystyle
-\frac{2}{n^2}\frac{\lambda_a^2 \lambda_b^2}{(\lambda_a-\lambda_b)^2}+O(n^{-3}),
&\text{if $a\ne b$.}
\end{cases}
\\
&
\nonumber\\
&E[(\hat\lambda_a-\lambda_a)(\hat\lambda_b-\lambda_b)]  (\triangleq V^{ab}(\hat\lam)) \nonumber\\
\label{asymp_var_result_mle}
&=
\begin{cases}
\displaystyle
\frac{2}{n}\lambda_a^2+O(n^{-5/2}), & \text{if $a=b$,}\\
&\\
\displaystyle
O(n^{-5/2}),
&\text{if $a\ne b$.}
\end{cases}
\end{align}
\end{prop}
This result says that $\hat{\lam}$ is the second-order efficient (among the bias-corrected first-order efficient estimators), but the bias-corrected sample eigenvalues are not. The difference in the asymptotic performance between the two estimators is due to the fact $\el^*$ do not use the prior information $\Gam=\Gam^0$, while $\hat\lam$ does. In contrast to $\el^*$, which does not use $\hh$,  $\hat\lam$ incorporates the information of  $\hh$ with the aid of the prior knowledge $\Gam=\Gam^0.$
In fact, as we will see in the next subsection, the difference between \eqref{asymp_var_result_sample_eigen} and \eqref{asymp_var_result_mle} is closely related to the asymptotic information loss caused by discarding  $\hh$.
\subsection{Asymptotic Information Loss}
\label{subsec:info loss}
In this subsection,  we consider the asymptotic information loss caused by ignoring $\hh$ for the estimation of $\lam$. Information loss matrix $(\Delta g_{ab}(\el)),\ 1\leq a, b\leq p$ at a fixed point $\s=(\lam, \Gam)$ is given by
$$
\Delta g_{ab}(\el)\triangleq E[g_{ab}(\bm{S} | \el)]=g_{ab}(\bm{S})-g_{ab}(\el),
$$
where $g_{ab}(\bm{S}),g_{ab}(\el), g_{ab}(\bm{S} | \el)$ are  the components of the metrics w.r.t. $\partial_a$ and $\partial_b$ based on respectively the distributions $\bm{S}$, $\el$ and the conditional distribution of $\bm{S}$ given $\el$, all of which are measured at the point $\s=(\lam, \Gam)$.

Amari \cite{Amari1} found that the asymptotic information loss can be expressed in terms of the metric and the embedding curvatures;
\begin{align}
\label{information_loss}
\Delta g_{ab}(\el)=&\;n \sum_{s<t, u<v}g_{a(s,t)}g_{b_(u,v)}g^{(s,t)(u,v)} \nonumber\\
&+\sum_{c,d,s<t,u<v}\stackrel{e}{H}_{ac(s,t)}\;\stackrel{e}{H}_{bd(u,v)}g^{cd}\;g^{(s,t)(u,v)}\nonumber\\
&+(1/2)\sum_{s<t,u<v,o<p,q<r}\stackrel{m}{H}_{(s,t)(u,v)a}\;\stackrel{m}{H}_{(o,p)(q,r)b}\;g^{(s,t)(o,p)}\;g^{(u,v)(q,r)}\nonumber\\
&+O(n^{-1}).
\end{align}
Straightforward calculation leads us to the following result:
\begin{prop}
\label{info_loss_result}
$$
\Delta g_{ab}(\el)=B_{ab}+O(n^{-1}),
$$
where 
$$
B_{ab}=
\begin{cases}
\displaystyle
\frac{1}{2\lambda_a^{2}}\sum_{t\ne a} \frac{\lambda_t^2}{(\lambda_t-\lambda_a)^2},
 &\text{if $a=b$},\\
\displaystyle
-\frac{1}{2(\lambda_a-\lambda_b)^2}, &\text{if $a\ne b$}.
\end{cases}
$$
\end{prop}

$B_{ab}$ at the point $(\lam, \Gam)$ depends only on $\lam.$ When the information loss of a statistic has the order $O(n^{-q+1})$, we call the statistic is the $q$th order sufficient. Consequently the statistic $\el$ is the first order sufficient, but not the second order sufficient. 

$B_{ab}$, the information loss in the second order term ($O(1)$) could be quite large when the population eigenvalues are close to each other. Note that the information carried by $\el$ is given by the formula;
\begin{align*}
g_{ab}(\el)&=g_{ab}(\bm{S})-\Delta g_{ab}(\el)\\
&=n g_{ab}(\xx)-\Delta g_{ab}(\el)\\
&=(n/2) \lambda_a^{-2}\delta(a=b)-\Delta g_{ab}(\el).
\end{align*}
Since $(g_{ab}(\el))$ is positive definite, $\text{diag}(n2^{-1}\lambda_1^{-2}, \ldots, n2^{-1}\lambda_p^{-2}) > (\Delta g_{ab})$. This holds true even in the neighborhood of a point $\lambda_1= \cdots =\lambda_p$ where $B_{ab}$ diverges. This indicates that the term of  order $O(n^{-1})$ in $\Delta g_{ab}(\el)$ is also unbounded in such a neighborhood. Hence the expansion of the information loss with respect to $n$ is not useful when the population eigenvalues are close to each other. 

Except for the case where the population eigenvalues are close to each other, Proposition \ref{info_loss_result} tells us approximately how much information is lost by ignoring the sample eigenvectors for the inference on the population eigenvalues. If we contract $\Delta g_{ab}$, then we could get a scalar measure on the information loss;
$$
IL\triangleq \sum_{a,b} g^{ab} \Delta g_{ab} = \sum_a  2 \lambda_a^2 B_{aa}+O(n^{-1}) = \sum_{a<b} \frac{\lambda_a^2+\lambda_b^2}{(\lambda_a-\lambda_b)^2}+O(n^{-1})
$$
Asymptotic information loss is closely related to the asymptotic variance of the two estimators $\el^*$ and $\hat\lambda$ in the previous subsection. Actually if we contract the asymptotic performance difference between the two estimators $V^{ab}(\el^*)-V^{ab}(\hat\lambda)$, then 
it equals $n^{-2}IL$, that is, 
\begin{align}
\label{IL&V}
&\sum_{a,b} (V^{ab}(\el^*)-V^{ab}(\hat\lam)) g_{ab}\nonumber \\
&=2^{-1}E[\sum_a (l^*_a/\lambda_a-1)^2]-2^{-1}E[\sum_a ({\hat\lambda}_a/\lambda_a-1)^2]\nonumber \\
&=n^{-2}\sum_{a<b} \frac{\lambda_a^2+\lambda_b^2}{(\lambda_a-\lambda_b)^2}+O(n^{-3}) =n^{-2} IL.
\end{align}
As a numerical example, we made a simulation for the case $p=2,\ n=20$. Taking the relationship \eqref{IL&V} into account, we could measure an information loss as the normalized quadratic risk difference between $\el^*$ and $\hat{\lam}.$ We randomly generated a two-dimensional normal vector under the following conditions, $\s=\text{diag}(1.0,\ c),\ c=0.2, 0.4, 0.6, 0.8, 1.0$. We made $10^8$ times repetition and took the average for each condition. 
The Table \ref{el_vs_hat_lam} shows the result. (Note: 1)The risk of $\hat\lam$ theoretically equals 0.4. 2)The simulated risk of $\el^*$ is quite unstable  as its large s.d. shows.) We notice that information loss is not negligible. The risk of $\el^*$ is larger than that of $\hat\lam$ by 24--111 \%. The risk difference is quite large especially when the population eigenvalues are close to each other.
\begin{table}
\caption{Simulate risk of $\el^*$ when $p=2$ as $c$ varies}
\label{el_vs_hat_lam}
\begin{center}
\begin{tabular}{l|r|r|r|r|r}
$c$ : Second Eigenvalue  & 1.0 & 0.8 & 0.6 & 0.4 & 0.2 \\
\hline
Simulated Risk of $\el^*$ & 0.85 &0.83 &0.70 &0.60 &0.50 \\
\hline
Standard Deviation & 0.24 & 0.48 &0.15 & 0.09 &0.22 \\
\hline
100$\times$(Risk Difference/Risk of $\hat\lam$) &111 &107 & 75 &49 & 24 \\
\end{tabular}
\end{center}
\end{table}
\section{Estimation of $\lam$ when $\Gam$ is unknown}
\label{section:estimation of lambda when Gamma is unknown}
In this section, we consider the more practical case where $\Gam$ is unknown. The derivation of a new estimator for this case will be done in view of the modification of the bias of $\bar\el$. Actually almost all the literature on the estimation of $\lam$ we mentioned in Section \ref{section:int}  modify the bias of $\bar\el$ by so called "shrinkage" method, that is, decreasing the dispersion of $\bar\el$. Though the concrete methods of shrinkage differ for each estimator, they are proposed mainly from analytical motivations. Here we consider another shrinkage estimator from a geometrical point of view.  

\begin{figure}[t]
\centering
\includegraphics[width=10cm,bb=9 9 200 200]{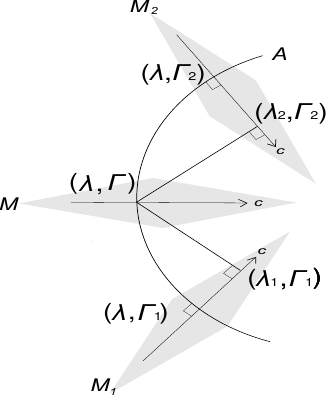}
\caption{ The shrinkage effect of the projection $(\lam, \Gam)$ onto $\mathcal{M}_i,\ i=1,2$}
\label{fig5}
\end{figure}

Suppose that we have a sample $\bar{\bm{S}}\triangleq n^{-1}\bm{S}$ which takes the point $(\lam, \Gam)$ in $\mathcal{S}$, that is, $\lam=\bar\el, \Gam=\bm{H}$. (See Figure \ref{fig5}.) Take the orthogonal projection of this point onto the submanifold $\mathcal{M}(\Gam_i)\triangleq \mathcal{M}_i (i=1,2)$, where the projected point $(\lam_i, \Gam_i)$ is given by $\lam_i=( (\Gam_i^t\bar{\bm{S}}\Gam_i)_{11},\ldots,(\Gam_i^t\bar{\bm{S}}\Gam_i)_{pp} )$. As we mentioned in Section \ref{section:embedding_curvature}, $(\lam_i, \Gam_i)$ is the minimum distance point on $\mathcal{M}_{i}$ from $(\lam, \Gam)$ with respect to Kullback-Leibler divergence. It is clearly understood that this projection has the shrinkage effect.  If we have an appropriate probability measure of $\Gam$ on the group of $p$-dimensional orthogonal matrices $\mathcal{O}(p)$, the expectation of $(\Gam^t\bar{\bm{S}}\Gam)_{ii}, i=1,\ldots,p$ for that measure would give birth to a natural shrinkage estimator.  

We choose the conditional distribution of $\hh$ when $\el$ is given for the probability measure on $\mathcal{O}(p)$. Since $\bm{S}=\hh \bm{L} \hh^t$ is distributed as Wishart matrix $W_p (n, \s)$, its density w.r.t. the uniform probability $d\mu(\hh)$ on $\mathcal{O}(p)$ equals 
\begin{equation}
f(\hh | \el\,;\s)= K(\el\,;\s)^{-1} \exp\Bigl(-(1/2)\tr \hh\bm{L} \hh^t \si\Bigr), 
\end{equation}
where normalizing constant $K(\el\,;\s)$ is given by
$$
K(\el\,;\s)=\int_{\mathcal{O}(p)} \exp\Bigl(-(1/2)\tr \hh\bm{L} \hh^t \si\Bigr) d\mu(\hh).
$$
This conditional distribution depends on $\s$. If we substitute $\s$ with an estimator $\hat{\s}(\bar{\bm{S}})$, 
it gives a distribution on $\mathcal{O}(p)$, whose density with respect to $d\mu(\Gam)$ is given by
\begin{equation}
\label{density_Gam}
f(\Gam | \el\,;\hat{\s})= K(\el\,;\bar{\bm{S}})^{-1} \exp\Bigl(-(1/2)\tr \Gam\bm{L} \Gam^t \hat{\s}^{-1}\Bigr), 
\end{equation}
where
$$
K(\el\,;\hat{\s})=\int_{\mathcal{O}(p)} \exp\Bigl(-(1/2)\tr \Gam\bm{L} \Gam^t\hat{\s}^{-1} \Bigr) d\mu(\Gam).
$$
Take the expectation of $(\Gam^t \bar{\bm{S}} \Gam)_{ii}$ w.r.t. the density \eqref{density_Gam}, then we have
\begin{equation}
\hat{\lambda}_i^{*}\triangleq K(\el\,;\bar{\bm{S}})^{-1} \int_{\mathcal{O}(p)} (\Gam^t \bar{\bm{S}} \Gam)_{ii} \exp\Bigl(-(1/2)\tr \Gam\bm{L} \Gam^t\hat{\s}^{-1} \Bigr) d\mu(\Gam),\quad i=1,\ldots,p.
\end{equation}
We propose $\hat{\lam}^*\triangleq (\hat{\lambda}_1^{*},\ldots,\hat{\lambda}_p^{*})$ as a new estimator of $\lam$. 

If $\hat\s$ is given by an orthogonally invariant estimator \eqref{orthogonally_inv_est},
$\hat{\lambda}_i^{*}$ can be more specifically described. Let $\bar{\bm{L}}$ denote $\diag(\bar\el)$. 
Because of the invariance of $d\mu$, it turns out that
\begin{eqnarray}
\label{new_est}
\hat{\lambda}_i^{*}&=& K(\el)^{-1} \int_{\mathcal{O}(p)} (\Gam^t \hh\bar{\bm{L}}\hh^t \Gam)_{ii} \exp\Bigl(-(1/2)\tr\bm{L} \Gam^t\hh\bm{\Phi}^{-1}\hh^t  \Gam\Bigr) d\mu(\Gam)\nonumber \\
&=&K(\el)^{-1} \int_{\mathcal{O}(p)} (\Gam^t \bar{\bm{L}} \Gam)_{ii} \exp\Bigl(-(1/2)\tr\bm{L} \Gam^t \bm{\Phi}^{-1}  \Gam\Bigr) d\mu(\Gam),
\end{eqnarray}
where
\begin{equation}
\label{normalizing_constant}
K(\el)\triangleq \int_{\mathcal{O}(p)}  \exp\Bigl(-(1/2)\tr\bm{L} \Gam^t \bm{\Phi}^{-1}  \Gam\Bigr) d\mu(\Gam).
\end{equation}
\begin{figure}[t]
\centering
\includegraphics[width=10cm,bb=9 9 458 434]{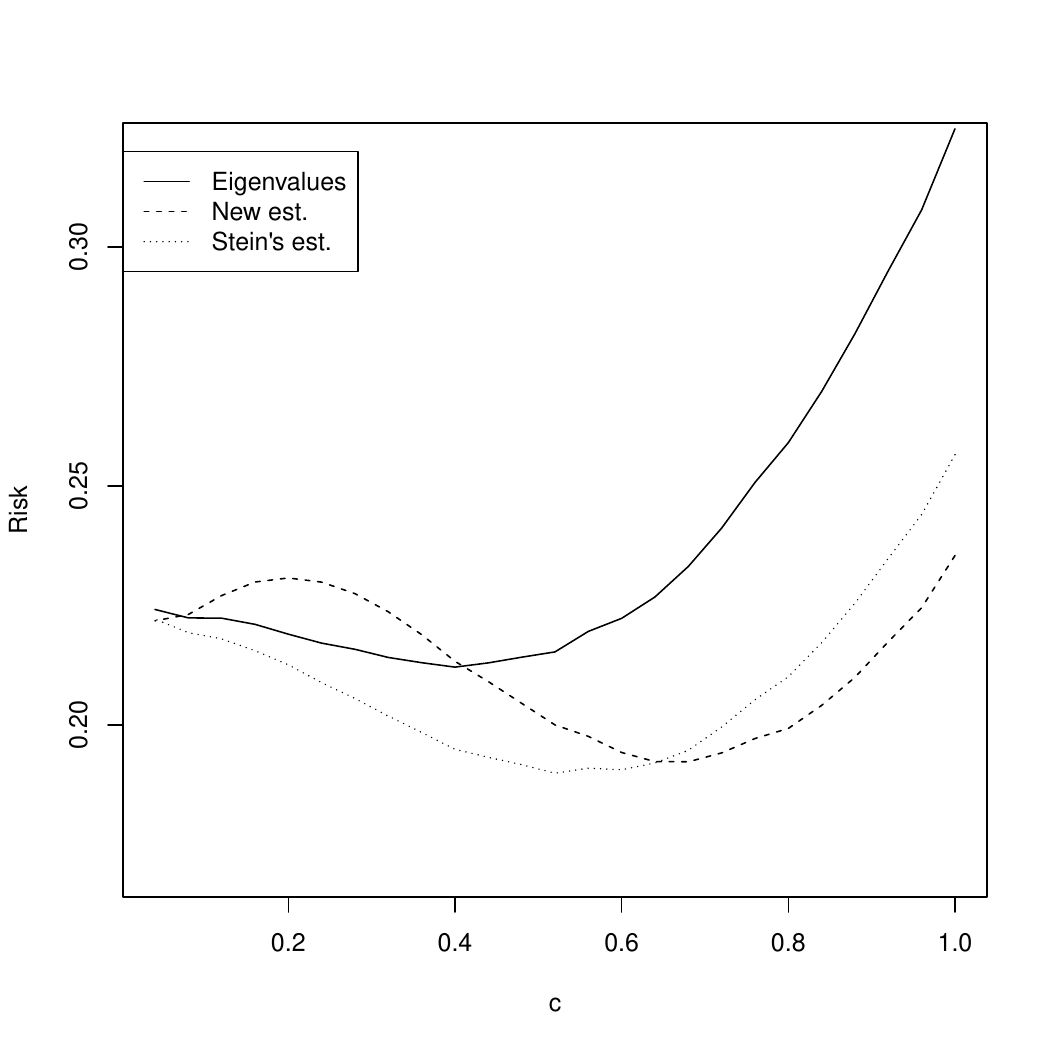}
\caption{ Risks of the three estimators as $c$ changes}
\label{fig6}
\end{figure}
The analytic evaluation of this estimator's performance seems difficult even for the large sample case. Instead we show the numerical result comparing $\bar{\el}$,  $\hat{\lam}^*$ and Stein's estimator \eqref{stein's_est}. Our new estimator $\hat{\lambda}_i^{*}$ is also equipped with the same $\phi$'s in \eqref{stein's_est}.
We simulated the risks of three estimators for the case $p=2,\ n=10$ w.r.t. K-L loss, which is given by 
$$
\sum_{i=1}^p \hat{\lambda}_i\lambda_i^{-1}-\sum_{i=1}^p \log(\hat{\lambda}_i\lambda_i^{-1})-p,
$$
where $\hat{\lambda}_i=\bar{l}_i,\ \hat{\lambda}^*_i,\ \phi_i,\ i=1,\ldots,p$.
Since all the estimators are functions of $\el$ and scale invariant, it is enough to measure the risks for  $\s=\diag(1, c),\ 0<c\leq 1$. We varied $c$ from 0.04 to 1.00 by the increment 0.04, and for each $c$ we repeated the risk evaluation $10^5$ times and took the average. For the integral calculation of \eqref{new_est} and \eqref{normalizing_constant}, we picked up 50 points from $\mathcal{O}(2)$ in an equidistant manner. Figure \ref{fig6} shows the result. The new estimator performs better compared to $\bar\el$, especially $\lam$ are close to each other, though it seems that $\hat{\lam}^*$ does not dominate $\bar\el$ as Stein's estimator does. Unfortunately we do not have any theoretical explanation of the risk behavior of the new estimator. We could only guess that the shrinkage effect works well when $c$ is close to one, while its effect is too strong elsewhere. We also simulated the risk of the new estimator equipped with M.L.E. instead of Stein's estimator. Since its performance  is almost the same as the above new estimator, we skip the result.
\section{Remark}
\begin{enumerate}
\item We treated the estimation problem of the eigenvalues $\lam$ in the latter half of the paper. The estimation on the eigenvectors $\Gam$ seems rather untouched in the classical situation $n\geq p$. Corollary \ref{efron_C} on the statistical curvatures of $\mathcal{A}$ or \eqref{asympto_bias} on the asymptotic bias tells us that the point where $\lam$ has some multiplicity is a statistically singular point. Around these points, inference on $\Gam$ are considered to need subtle treatment. Especially the eigenvectors are not well identified around the multiplicity point, hence the information contained in $\hh$ vanishes there (see $g_{(s,t)(u,v)}$ in Proposition \ref{metric_spectra}). This indicates that the inference using only $\hh$ is not appropriate.

\item We proposed a new estimator for $\lam$ in Section \ref{section:estimation of lambda when Gamma is unknown} . However this belongs to the same category as most estimators in the past literature in that it uses sample eigenvalues $\lam$ only. It is still unclear how we can use the sample eigenvalues $\hh$ for the inference of $\lam$.

\end{enumerate}
\section{Appendix}
\subsection{Proof of Proposition \ref{metric_spectra}}
\label{proof_prop1}
As a base for the vector space of real symmetric matrices, we consider $\bm{E}_{ij}\ (1\leq i \leq j \leq p)$ which is a $p\times p$ matrix defined by
$$
\bm{E}_{ij}=
\begin{cases}
\bm{I}_{ii},  &\text{if $i=j$},\\
\bm{I}_{ij}+\bm{I}_{ji},  &\text{if $i < j$},
\end{cases}
$$
where $\bm{I}_{ij}\ (1\leq i, j \leq p)$ is the $p \times p$ matrix whose $(i, j)$ element equals one, and all the other elements are zero. The one to one correspondence
$$
\partial^{(i,j)}\triangleq \frac{\partial \ \ }{\partial \sigma_{ij}} \longleftrightarrow \bm{E}_{ij}, \qquad 1\leq i\leq j \leq p,
$$
gives the component expression of \eqref{metric_simple_form}
$$
\langle \partial^{(i,j)},  \partial^{(k,l)} \rangle=\frac{1}{2} \text{tr}\Bigl(\si \bm{E}_{ij} \si \bm{E}_{kl}\Bigr),\qquad 1\leq i \leq j \leq p,\ 1\leq k \leq l \leq p.
$$
Since
\begin{align}
\partial_a&\triangleq \frac{\partial\ \ }{\partial \lambda_a}=\sum_{i \leq j}\frac{\partial \sigma_{ij}}{\partial \lambda_a}\frac{\partial\ \ }{\partial \sigma_{ij}}
=\sum_{i \leq j}\frac{\partial \sigma_{ij}}{\partial \lambda_a}\partial^{(i,j)}
\quad 1\leq a \leq p,\\
\partial_{(s,t)}&\triangleq \frac{\partial\ \ }{\partial u_{st}}=\sum_{i\leq j}\frac{\partial \sigma_{ij}}{\partial u_{st}}\frac{\partial\ \ }{\partial \sigma_{ij}}
=\sum_{i\leq j}\frac{\partial \sigma_{ij}}{\partial u_{st}}\partial^{(i,j)}
\quad 1\leq s<t \leq p,
\end{align}
we have the following relations
\begin{align}
\label{g_a_b}
g_{ab}&=\frac{1}{2}\tr\Bigl\{
\si\Bigl(\sum_{i\leq j} \frac{\partial \sigma_{ij}}{\partial \lambda_a} E_{ij}\Bigr) \si
\Bigl(\sum_{k\leq l} \frac{\partial \sigma_{kl}}{\partial \lambda_b} E_{kl}\Bigr)
\Bigr\}, \\
\label{g_a_(s,t)}
g_{a(s,t)}&=\frac{1}{2}\tr\Bigl\{
\si\Bigl(\sum_{i\leq j} \frac{\partial \sigma_{ij}}{\partial \lambda_a} E_{ij}\Bigr) \si
\Bigl(\sum_{k\leq l} \frac{\partial \sigma_{kl}}{\partial u_{st}} E_{kl}\Bigr)
\Bigr\},\\
\label{g_(s,t)_(u,v)}
g_{(s,t)(u,v)}&=\frac{1}{2}\tr\Bigl\{
\si\Bigl(\sum_{i\leq j} \frac{\partial \sigma_{ij}}{\partial u_{st}} E_{ij}\Bigr) \si
\Bigl(\sum_{k\leq l} \frac{\partial \sigma_{kl}}{\partial u_{uv}} E_{kl}\Bigr)
\Bigr\},
\end{align}
where $1\leq a, b \leq p$, $1\leq s <t \leq p,\  1 \leq u < v \leq p.$

For the first order derivative at $\uu=\bm{0}$, we only have to consider $\s$ up to the term to the first power w.r.t. $\uu$, hence we put $\s(\lam, \uu)$ as
\begin{align}
\s(\lam, \uu)&=\Gam(\bm{I}_p+\bm{U})\bm{\Lambda}(\bm{I}_p+\bm{U})^t\Gam^t +O(||\uu||^2)\nonumber \\
&=\Gam\bm{\Lambda}\Gam^t+\Gam\bm{\Lambda}\bm{U}^t \Gam^t +\Gam \bm{U}\bm{\Lambda}\Gam^t+O(||\uu||^2).
\end{align}
Therefore we have
$$
\sigma_{ij}=\sum_{k} \gamma_{ik}\gamma_{jk}\lambda_k
+\sum_{k,l}\gamma_{ik}\gamma_{jl}\Bigl(\lambda_{k}u_{lk}+\lambda_{l}u_{kl}\Bigr)+O(||\uu||^2), \quad 1\leq i\leq  j \leq p,
$$
where $u_{ii}\triangleq 0\ (1\leq i \leq p),\ u_{ij}\triangleq -u_{ji}\ (1\leq j<i \leq p),$
which leads to
\begin{equation}
\label{d_sigma/d_lambda}
\frac{\partial \sigma_{ij}}{\partial\lambda_a}\bigg|_{\uu=\bm{0}}=\gamma_{ia}\gamma_{ja},
\end{equation}
and
\begin{equation}
\label{d_sigma/d_u}
\frac{\partial \sigma_{ij}}{\partial u_{st}}\bigg|_{\uu=\bm{0}}
=\lambda_t\gamma_{it}\gamma_{js}-\lambda_s\gamma_{is}\gamma_{jt}
+\lambda_t\gamma_{is}\gamma_{jt}-\lambda_s\gamma_{it}\gamma_{js}.
\end{equation}
From \eqref{d_sigma/d_lambda} and \eqref{d_sigma/d_u}, we have the following results on tangent vectors;
\begin{equation}
\label{tangent_vector_eigenvalues}
\sum_{i\leq j} \frac{\partial \sigma_{ij}}{\partial \lambda_a} E_{ij}=\sum_{i \leq j}\gamma_{ia}\gamma_{ja} \bm{E}_{ij}= \bm{\gamma}_a\bm{\gamma}_a^t,
\end{equation}
where $\bm{\gamma}_a$ is the $a$th column of $\Gam$, and
\begin{equation}
\label{tangent_vector_eigenvectors}
\sum_{i\leq j} \frac{\partial \sigma_{ij}}{\partial u_{st}} E_{ij}=
\lambda_t\bm{\gamma}_t\bm{\gamma}_s^t-\lambda_s\bm{\gamma}_s\bm{\gamma}_t^t
+\lambda_t\bm{\gamma}_s\bm{\gamma}_t^t-\lambda_s\bm{\gamma}_t\bm{\gamma}_s^t.
\end{equation}
If we substitute \eqref{tangent_vector_eigenvalues} and \eqref{tangent_vector_eigenvectors} into \eqref{g_a_b}, \eqref{g_a_(s,t)} and \eqref{g_(s,t)_(u,v)}, we get the results as follows;
\begin{align*}
2g_{ab}&=\tr\Bigl(\si  \bm{\gamma}_a\bm{\gamma}_a^t \si  \bm{\gamma}_b\bm{\gamma}_b^t \Bigr)\\
&=\tr\Bigl\{\Bigl(\bm{\gamma}_b^t \si \bm{\gamma}_a\Bigr)\Bigl\{\Bigl(\bm{\gamma}_a^t \si \bm{\gamma}_b\Bigr)\Bigr\}\\
&=\tr\Bigl(\lambda_a^{-1} \delta(a=b) \lambda_b^{-1} \delta(a=b)\Bigr)\\
&=\lambda_a^{-2}\delta(a=b),
\end{align*}
\begin{align*}
2g_{a(s,t)}&=\tr\Bigl\{\si \bm{\gamma}_a\bm{\gamma}_a^t \si 
\Bigl(
\lambda_t\bm{\gamma}_t\bm{\gamma}_s^t-\lambda_s\bm{\gamma}_s\bm{\gamma}_t^t
+\lambda_t\bm{\gamma}_s\bm{\gamma}_t^t-\lambda_s\bm{\gamma}_t\bm{\gamma}_s^t
\Bigr)\Bigr\}\\
&=\lambda_t \lambda_a^{-2} \delta(a=s=t)-\lambda_s\lambda_a^{-2}\delta(a=s=t)\\ \Bigr.
&\qquad\Bigl.+\lambda_t\lambda_a^{-2}\delta(a=s=t)-\lambda_s\lambda_a^{-2}\delta(a=s=t)\\
&=0,\\
2g_{(s,t)(u,v)}&=\tr\Bigl\{
\si\Bigl(
\lambda_t\bm{\gamma}_t\bm{\gamma}_s^t-\lambda_s\bm{\gamma}_s\bm{\gamma}_t^t
+\lambda_t\bm{\gamma}_s\bm{\gamma}_t^t-\lambda_s\bm{\gamma}_t\bm{\gamma}_s^t
\Bigr)
\Bigl.\\
&\qquad\times\Bigl.\si \Bigl(\lambda_v\bm{\gamma}_v\bm{\gamma}_u^t-\lambda_u\bm{\gamma}_u\bm{\gamma}_v^t
+\lambda_v\bm{\gamma}_u\bm{\gamma}_v^t-\lambda_u\bm{\gamma}_v\bm{\gamma}_u^t
\Bigr)
\Bigr\}\\
&=\lambda_t\lambda_v\lambda_t^{-1}\delta(u=t)\lambda_s^{-1}\delta(s=v)
-\lambda_t\lambda_u\lambda_t^{-1}\delta(v=t)\lambda_s^{-1}\delta(s=u)\\
&\qquad+\lambda_t\lambda_v\lambda_t^{-1}\delta(v=t)\lambda_s^{-1}\delta(s=u)
-\lambda_t\lambda_u\lambda_t^{-1}\delta(u=t)\lambda_s^{-1}\delta(s=v)\\
&\qquad-\lambda_s\lambda_v\lambda_s^{-1}\delta(u=s)\lambda_t^{-1}\delta(t=v)
+\lambda_s\lambda_u\lambda_s^{-1}\delta(s=v)\lambda_t^{-1}\delta(t=u)\\
&\qquad-\lambda_s\lambda_v\lambda_s^{-1}\delta(s=v)\lambda_t^{-1}\delta(t=u)
+\lambda_s\lambda_u\lambda_s^{-1}\delta(u=s)\lambda_t^{-1}\delta(t=v)\\
&\qquad+\lambda_t\lambda_v\lambda_s^{-1}\delta(u=s)\lambda_t^{-1}\delta(t=v)
-\lambda_t\lambda_u\lambda_s^{-1}\delta(v=s)\lambda_t^{-1}\delta(t=u)\\
&\qquad+\lambda_t\lambda_v\lambda_s^{-1}\delta(v=s)\lambda_t^{-1}\delta(t=u)
-\lambda_t\lambda_u\lambda_s^{-1}\delta(u=s)\lambda_t^{-1}\delta(t=v)\\
&\qquad-\lambda_s\lambda_v\lambda_t^{-1}\delta(u=t)\lambda_s^{-1}\delta(s=v)
+\lambda_s\lambda_u\lambda_t^{-1}\delta(v=t)\lambda_s^{-1}\delta(u=s)\\
&\qquad-\lambda_s\lambda_v\lambda_t^{-1}\delta(t=v)\lambda_s^{-1}\delta(u=s)
+\lambda_s\lambda_u\lambda_t^{-1}\delta(u=t)\lambda_s^{-1}\delta(s=v)\\
&=(-1+\lambda_t\lambda_s^{-1}-1+\lambda_s\lambda_t^{-1}+\lambda_t\lambda_s^{-1}-1+\lambda_s\lambda_t^{-1}-1)\delta(s=u,t=v)\\
&=2(\lambda_s^{-1}(\lambda_t-\lambda_s)+\lambda_t^{-1}(\lambda_s-\lambda_t))\delta(s=u,t=v)\\
&=2(\lambda_t-\lambda_s)(\lambda_s^{-1}-\lambda_t^{-1})\delta(s=u,t=v)\\
&=2(\lambda_t-\lambda_s)^2(\lambda_s\lambda_t)^{-1}\delta(s=u,t=v).
\end{align*}
\subsection{Proof of Proposition \ref{embedding_c}}
Note that $\si=\Gam \bm{\Lambda}^{-1} \Gam^t$, hence
$$
\theta^{ij}=
\begin{cases}
\displaystyle
-\sum_k \gamma_{ik} \gamma_{jk} \lambda_k^{-1} &\text{ if $i < j$,} \\
\displaystyle
-2^{-1} \sum_k \gamma_{ik}^2 \lambda_k^{-1} &\text{ if $i=j$.}
\end{cases}
$$
This means $\mathcal{M}$ is an affine subspace of $\mathcal{S}$ w.r.t. an $\tteta$, which is an affine coordinate system of $\mathcal{S}$ with $e$-connection. Consequently $\mathcal{M}$ is $e$-flat, i.e. $\stackrel{e}{H}_{ab(s,t)}=0$. $\stackrel{m}{H}_{ab(s,t)}=0$ is similarly proved. See Theorem 1.1 in Amari and Nagaoka \cite{Amari&Nagaoka}. 

Now we consider $\stackrel{m}{H}_{(s,t)(u,v)a}$. Using (4.14) in Amari \cite{Amari2}, it is calculated as
\begin{align*}
\label{embedding_m-curvature_M_calculation}
\stackrel{m}{H}_{(s,t)(u,v)a}&=\sum_{i\leq j} \frac{\partial^2 \sigma_{ij}}{\partial u_{st}\partial u_{uv}}\bigg|_{\uu=\bm{0}}\frac{\partial \theta^{ij}}{\partial \lambda_a}\bigg|_{\uu=\bm{0}}\\
&=-2^{-1}\sum_{1\leq i, j \leq p} \frac{\partial^2 \sigma_{ij}}{\partial u_{st}\partial u_{uv}}\bigg|_{\uu=\bm{0}}\frac{\partial \sigma^{ij}}{\partial \lambda_a}\bigg|_{\uu=\bm{0}}\\
&=-2^{-1}\tr (\bm{A}\bm{B}),
\end{align*}
where $p\times p$ matrices $\bm{A}$, $\bm{B}$ are given by
$$
(\bm{A})_{ij}\triangleq\frac{\partial^2 \sigma_{ij}}{\partial u_{st}\partial u_{uv}}\bigg|_{\uu=\bm{0}},\
(\bm{B})_{ij}\triangleq\frac{\partial \sigma^{ij}}{\partial \lambda_a}\bigg|_{\uu=\bm{0}},\quad 1\leq i, j \leq p.
$$
In order to calculate $\bm{A}$, we only have to consider $\s$ up to the terms powered by two w.r.t. $\uu$;
\begin{align*}
\s&=\Gam\Bigl(\bm{I}_p+\bm{U}+2^{-1}\bm{U}^2\Bigr)\bm{\Lambda}\Bigl(\bm{I}_p+\bm{U}+2^{-1}\bm{U}^2\Bigr)^t\Gam^t
+O(||\uu||^3)\\
&=\Gam \bm{\Lambda}\Gam^t +\Gam(\bm{U}\bm{\Lambda}+\bm{\Lambda}\bm{U}^t)\Gam^t
+2^{-1}\Gam(\bm{U}^2\bm{\Lambda}+\bm{\Lambda}(\bm{U}^2)^t)\Gam^t+\Gam\bm{U}\bm{\Lambda}\bm{U}^t\Gam^t\\
&\qquad+O(||\uu||^3).
\end{align*}
Therefore $\sigma_{ij}$ is expressed as
\begin{equation}
\label{sigma^2_expression_upto_second_order}
\sigma_{ij}=2^{-1}\sum_{k,l}\gamma_{ik}\gamma_{jl}\Bigl((\bm{U}^2\bm{\Lambda}+\bm{\Lambda}(\bm{U}^2)^t)_{kl}+2(\bm{U}\bm{\Lambda}\bm{U}^t)_{kl}\Bigr)+\bm{R}_{ij}+O(||\uu||^3),
\end{equation}
where $\bm{R}=\Gam \bm{\Lambda}\Gam^t +\Gam(\bm{U}\bm{\Lambda}+\bm{\Lambda}\bm{U}^t)\Gam^t.$
Since
\begin{align*}
(\bm{U}^2\bm{\Lambda}+\bm{\Lambda}(\bm{U}^2)^t)_{kl}&=(\bm{U}^2\bm{\Lambda})_{kl}+(\bm{U}^2\bm{\Lambda})_{lk}\\
&=\sum_{b} u_{kb} u_{bl} \lambda_l + \sum_{b} u_{lb} u_{bk} \lambda_k,\\
2(\bm{U}\bm{\Lambda}\bm{U}^t)_{kl}&=2\sum_b u_{kb}u_{lb}\lambda_b,
\end{align*}
\eqref{sigma^2_expression_upto_second_order} truns out to be
\begin{equation}
\label{sigma^2_expression_upto_second_order_ver2}
\sigma_{ij}=2^{-1}\sum_{k,l,b} \gamma_{ik}\gamma_{jl}(u_{kb}u_{bl}\lambda_l+u_{lb}u_{bk}\lambda_k+2u_{kb}u_{lb}\lambda_b)+\bm{R}_{ij}+O(||\uu||^3).
\end{equation}
From this we have
\begin{equation}
\frac{\partial^2 \sigma_{ij}}{\partial u_{st}\partial u_{uv}}\bigg|_{\uu=\bm{0}}\times 2=\sum_{k,l,b}(a^{(1)}_{ij}+a^{(1)}_{ji}+a^{(2)}_{ij}+a^{(2)}_{ji}+a^{(3)}_{ij}+a^{(4)}_{ij}),
\end{equation}
where
\begin{align*}
a^{(1)}_{ij}
&=\gamma_{is} \gamma_{jv} \lambda_v \delta\{(k,b)=(s,t), (b,l)=(u,v), (s,t)\ne(u,v)\}\\
&\quad -\gamma_{it} \gamma_{jv} \lambda_v \delta\{(k,b)=(t,s), (b,l)=(u,v), (s,t)\ne(u,v)\}\\
&\quad -\gamma_{is} \gamma_{ju} \lambda_u \delta\{(k,b)=(s,t), (b,l)=(v,u), (s,t)\ne(u,v)\}\\
&\quad +\gamma_{it} \gamma_{ju} \lambda_u \delta\{(k,b)=(t,s), (b,l)=(v,u), (s,t)\ne(u,v)\}\\
&\quad +\gamma_{iu} \gamma_{jt} \lambda_t \delta\{(k,b)=(u,v), (b,l)=(s,t), (s,t)\ne(u,v)\}\\
&\quad -\gamma_{iv} \gamma_{jt} \lambda_t \delta\{(k,b)=(v,u), (b,l)=(s,t), (s,t)\ne(u,v)\}\\
&\quad -\gamma_{iu} \gamma_{js} \lambda_s \delta\{(k,b)=(u,v), (b,l)=(t,s), (s,t)\ne(u,v)\}\\
&\quad +\gamma_{iv} \gamma_{js} \lambda_s \delta\{(k,b)=(v,u), (b,l)=(t,s), (s,t)\ne(u,v)\},
\end{align*}
\begin{align*}
a^{(2)}_{ij}
&=2\gamma_{is}\gamma_{jt}\lambda_t\delta\{(k,b)=(s,t),(b,l)=(s,t),(s,t)=(u,v)\}\\
&\quad +2\gamma_{it}\gamma_{js}\lambda_s\delta\{(k,b)=(t,s),(b,l)=(t,s),(s,t)=(u,v)\}\\
&\quad -2\gamma_{is}\gamma_{js}\lambda_s\delta\{(k,b)=(s,t),(b,l)=(t,s),(s,t)=(u,v)\}\\
&\quad -2\gamma_{it}\gamma_{jt}\lambda_t\delta\{(k,b)=(t,s),(b,l)=(s,t),(s,t)=(u,v)\}\\
&=-2\gamma_{is}\gamma_{js}\lambda_s\delta\{(k,b)=(s,t),(b,l)=(t,s),(s,t)=(u,v)\}\\
&\quad -2\gamma_{it}\gamma_{jt}\lambda_t\delta\{(k,b)=(t,s),(b,l)=(s,t),(s,t)=(u,v)\},
\end{align*}
\begin{align*}
a^{(3)}_{ij}
&=2\gamma_{is}\gamma_{ju}\lambda_t\delta\{(k,b)=(s,t),(l,b)=(u,v),(s,t)\ne (u,v)\}\\
&\quad -2\gamma_{it}\gamma_{ju}\lambda_s\delta\{(k,b)=(t,s),(l,b)=(u,v),(s,t)\ne (u,v)\}\\
&\quad -2\gamma_{is}\gamma_{jv}\lambda_t\delta\{(k,b)=(s,t),(l,b)=(v,u),(s,t)\ne (u,v)\}\\
&\quad +2\gamma_{it}\gamma_{jv}\lambda_s\delta\{(k,b)=(t,s),(l,b)=(v,u),(s,t)\ne (u,v)\}\\
&\quad +2\gamma_{iu}\gamma_{js}\lambda_t\delta\{(k,b)=(u,v),(l,b)=(s,t),(s,t)\ne (u,v)\}\\
&\quad -2\gamma_{iv}\gamma_{js}\lambda_t\delta\{(k,b)=(v,u),(l,b)=(s,t),(s,t)\ne (u,v)\}\\
&\quad -2\gamma_{iu}\gamma_{jt}\lambda_s\delta\{(k,b)=(u,v),(l,b)=(t,s),(s,t)\ne (u,v)\}\\
&\quad +2\gamma_{iv}\gamma_{jt}\lambda_s\delta\{(k,b)=(v,u),(l,b)=(t,s),(s,t)\ne (u,v)\},
\end{align*}
\begin{align*}
a^{(4)}_{ij}
&= 4\gamma_{is}\gamma_{js}\lambda_t\delta\{(k,b)=(l,b)=(s,t)=(u,v)\}\\
&\quad +4\gamma_{it}\gamma_{jt}\lambda_s\delta\{(k,b)=(l,b)=(t,s)=(v,u)\}\\
&\quad -4\gamma_{is}\gamma_{jt}\lambda_t\delta\{(k,b)=(s,t), (l,b)=(t,s), (s,t)=(u,v)\}\\
&\quad -4\gamma_{it}\gamma_{js}\lambda_s\delta\{(k,b)=(t,s), (l,b)=(s,t), (s,t)=(u,v)\}\\
&= 4\gamma_{is}\gamma_{js}\lambda_t\delta\{(k,b)=(l,b)=(s,t)=(u,v)\}\\
&\quad +4\gamma_{it}\gamma_{jt}\lambda_s\delta\{(k,b)=(l,b)=(t,s)=(v,u)\}.
\end{align*}
Furthermore we have
\begin{equation}
\label{form_A}
2\bm{A}=\bm{A}^{(1)}+(\bm{A}^{(1)})^t+\bm{A}^{(2)}+(\bm{A}^{(2)})^t+\bm{A}^{(3)}+\bm{A}^{(4)},
\end{equation}
where
\begin{align*}
\bm{A}^{(1)}&=
\bm{\gamma}_s\bm{\gamma}_v^t\lambda_v \delta(t=u)+\bm{\gamma}_t\bm{\gamma}_u^t\lambda_u \delta(s=v)\\
&\quad +\bm{\gamma}_u\bm{\gamma}_t^t\lambda_t \delta(s=v)+\bm{\gamma}_v\bm{\gamma}_s^t\lambda_s \delta(u=t)\\
&\quad -\bm{\gamma}_t\bm{\gamma}_v^t\lambda_v \delta(s=u,t\ne v)-\bm{\gamma}_s\bm{\gamma}_u^t\lambda_u \delta(t=v,s\ne u)\\
&\quad -\bm{\gamma}_v\bm{\gamma}_t^t\lambda_t \delta(u=s, t\ne v)-\bm{\gamma}_u\bm{\gamma}_s^t\lambda_s \delta(t=v, s\ne u),\\
\bm{A}^{(2)}&=
-2(\bm{\gamma}_s \bm{\gamma}_s^t \lambda_s \delta(s=u, t=v)+\bm{\gamma}_t \bm{\gamma}_t^t \lambda_t \delta(s=u, t=v)),\\
\bm{A}^{(3)}&=2\Bigl(\bm{\gamma}_s\bm{\gamma}_u^t\lambda_t \delta(t=v, s\ne u)+\bm{\gamma}_t\bm{\gamma}_v^t\lambda_s \delta(s=u, t\ne v)\\
&\qquad+\bm{\gamma}_u\bm{\gamma}_s^t\lambda_t\delta(v=t, s\ne u)+\bm{\gamma}_v\bm{\gamma}_t^t\lambda_s \delta(u=s, t\ne v)\Bigr)\\
&\qquad -2\Bigl(\bm{\gamma}_t\bm{\gamma}_u^t\lambda_s \delta(s=v)+\bm{\gamma}_s\bm{\gamma}_v^t\lambda_t \delta(t=u)+\bm{\gamma}_v\bm{\gamma}_s^t\lambda_t\delta(u=t)
+\bm{\gamma}_u\bm{\gamma}_t^t\lambda_s \delta(s=v)\Bigr),\\
\bm{A}^{(4)}&=4\Bigl(\bm{\gamma}_s\bm{\gamma}_s^t\lambda_t \delta(s=u, t=v)
+\bm{\gamma}_t\bm{\gamma}_t^t\lambda_s\delta(s=u, t=v)\Bigr).
\end{align*}
Since
$$
\frac{\partial \sigma^{ij}}{\partial \lambda_a}\Bigg|_{\uu=\bm{0}}=-\lambda_a^{-2}\gamma_{ia}\gamma_{ja},
$$
we have
\begin{equation}
\label{form_B}
\bm{B}=-\lambda_a^{-2} \bm{\gamma}_a \bm{\gamma}_a^t.
\end{equation}
From \eqref{form_A} and \eqref{form_B}, we have
\begin{align*}
\stackrel{m}{H}_{(s,t)(u,v)a}&=-4^{-1}\tr(2\bm{A}\bm{B})\\
&=-4^{-1}\tr\{(\bm{A}^{(1)}+(\bm{A}^{(1)})^t+\bm{A}^{(2)}+(\bm{A}^{(2)})^t+\bm{A}^{(3)}+\bm{A}^{(4)})\bm{B}\}\\
&=4^{-1}\lambda_a^{-2}\tr\{(\bm{A}^{(1)}+(\bm{A}^{(1)})^t+\bm{A}^{(2)}+(\bm{A}^{(2)})^t+\bm{A}^{(3)}+\bm{A}^{(4)})
\bm{\gamma}_a \bm{\gamma}_a^t\}.
\end{align*}
The following equalities hold;
\begin{align*}
\tr(\bm{A}^{(1)}\bm{\gamma}_a\bm{\gamma}_a^t)
&=\lambda_a\delta(s=v=a,t=u)+\lambda_a\delta(t=u=a,s=v)\\
&\quad+\lambda_a\delta(t=u=a,s=v)+\lambda_a\delta(s=v=a,t=u)\\
&\quad-\lambda_a\delta(t=v=a,s=u,t\ne v)-\lambda_a\delta(s=u=a,t=v,s\ne u)\\
&\quad-\lambda_a\delta(t=v=a,s=u, t\ne v)-\lambda_a\delta(s=u=a,t=v,s\ne u)\\
&=0. \\
\tr((\bm{A}^{(1)})^t\bm{\gamma}_a\bm{\gamma}_a^t)
&=0.\\
\tr(\bm{A}^{(2)}\bm{\gamma}_a\bm{\gamma}_a^t)&=-2(\lambda_a \delta(s=u=a, t=v)+\lambda_a\delta(s=u, t=v=a)).\\
\tr((\bm{A}^{(2)})^t\bm{\gamma}_a\bm{\gamma}_a^t)&=-2(\lambda_a \delta(s=u=a, t=v)+\lambda_a\delta(s=u, t=v=a)).\\
\tr(\bm{A}^{(3)}\bm{\gamma}_a\bm{\gamma}_a^t)
&=2\{\lambda_t\delta(s=u=a)\delta(t=v,s\ne  u)+\lambda_s\delta(t=v=a)\delta(s=u,t\ne v)\\
&\quad+\lambda_t\delta(s=u=a)\delta(t=v,s\ne  u)+\lambda_s\delta(t=v=a)\delta(s=u,t\ne v)\}\\
&\quad-2\{\lambda_s\delta(t=u=a)\delta(s=v)
+\lambda_t\delta(s=v=a)\delta(t=u)\\
&\quad+\lambda_t\delta(s=v=a)\delta(t=u)
+\lambda_s\delta(u=t=a)\delta(s=v)\}\\
&=0.\\
\tr(\bm{A}^{(4)}\bm{\gamma}_a\bm{\gamma}_a^t)&=4\lambda_t \delta(s=u=a,t=v)+4\lambda_s\delta(t=v=a,s=u).
\end{align*}
Consequently 
\begin{align*}
\stackrel{m}{H}_{(s,t)(u,v)a}
&=-\lambda_a^{-1}\delta(s=u=a,t=v)-\lambda_a^{-1}\delta(s=u,t=v=a)\\
&\quad +\lambda_a^{-2}\lambda_t \delta(s=u=a,t=v)+\lambda_a^{-2}\lambda_s\delta(t=v=a,s=u)\\
&=
\begin{cases}
\lambda_a^{-2}(\lambda_t-\lambda_a),  \text{ if $s=u=a$, $t=v$,}\\
\lambda_a^{-2}(\lambda_s-\lambda_a),  \text{ if $s=u$, $t=v=a$,}\\
0, \text{ otherwise.}
\end{cases}
\end{align*}
\subsection{Proof of Corollary \ref{efron_C}}
As we will see in the next subsection, 
\begin{align*}
&\sum_{s<t,u<v,o<p,q<r}\stackrel{m}{H}_{(s,t)(u,v)a}\;\stackrel{m}{H}_{(o,p)(q,r)b}\;g^{(s,t)(o,p)}\;g^{(u,v)(q,r)}\\
&=
\begin{cases}
\displaystyle
\frac{1}{\lambda_a^{2}}\sum_{t\ne a}\frac{\lambda_t^2}{(\lambda_t-\lambda_a)^2}, \text{ if $a=b$,}\\
\displaystyle
-\frac{1}{(\lambda_a-\lambda_a)^2}, \text{ if $a\ne b$.}
\end{cases}
\end{align*}
Combine this with Proposition \ref{metric_spectra}, we have 
\begin{align*}
\displaystyle
\gamma(\mathcal{A})&=2\sum_a \sum_{t\ne a}\frac{\lambda_t^2}{(\lambda_t-\lambda_a)^2}\\
\displaystyle
&=2\sum_{a<b}\frac{\lambda_a^2+\lambda_b^2}{(\lambda_a-\lambda_b)^2}.
\end{align*}
\subsection{Proof of Proposition \ref{asymp_var_result}}
We calculate each term in \eqref{asympto_variance}.  $g^{ab}=\delta(a=b)2\lambda_a^2$ from Proposition \ref{metric_spectra}.
Because of \eqref{e-embedding_c_M_upper_comp_vanish} and \eqref{conection_coeff_M_vanish}, 
$$
(\varGamma_M^{m})^{2ab}=(H_M^e)^{2ab}=0.
$$
For $\el^*$, $(H_{\hat{A}}^m)^{2ab}=(H_{A}^m)^{2ab}$. 
\begin{align}
\label{H_A^m_cal}
(H_{A}^m)^{2ab}&=\sum_{s<t, u<v, o<p, q<r}\mathcomposite[3]{\mathrel}{\hspace{-11mm} {\scriptsize \textit m}}{H^{a}_{(s,t)(u,v)}} \:\mathcomposite[1]{\mathrel}{\hspace{-11mm}{\scriptsize\textit m}}{H^{b}_{(o,p)(q,r)}} g^{(s,t)(o,p)}g^{(u,v)(q,r)}\nonumber\\
&=\sum_{t>a,p>b} \mathcomposite[3]{\mathrel}{\hspace{-11mm} {\scriptsize\textit m}}{H^{a}_{(a,t)(a,t)}} \:\mathcomposite[1]{\mathrel}{\hspace{-11mm} {\scriptsize\textit m}}{H^{b}_{(b,p)(b,p)}} (g^{(a,t)(b,p)})^2\nonumber\\
&+\sum_{t>a,p<b} \mathcomposite[3]{\mathrel}{\hspace{-11mm} {\scriptsize\textit m}}{H^{a}_{(a,t)(a,t)}} \:\mathcomposite[1]{\mathrel}{\hspace{-11mm} {\scriptsize\textit m}}{H^{b}_{(p,b)(p,b)}} (g^{(a,t)(p,b)})^2\nonumber\\
&+\sum_{t<a,p>b} \mathcomposite[3]{\mathrel}{\hspace{-11mm} {\scriptsize\textit m}}{H^{a}_{(t,a)(t,a)}} \:\mathcomposite[1]{\mathrel}{\hspace{-11mm} {\scriptsize\textit m}}{H^{b}_{(b,p)(b,p)}} (g^{(t,a)(b,p)})^2\nonumber\\
&+\sum_{t<a,p<b} \mathcomposite[3]{\mathrel}{\hspace{-11mm} {\scriptsize\textit m}}{H^{a}_{(t,a)(t,a)}} \:\mathcomposite[1]{\mathrel}{\hspace{-11mm} {\scriptsize\textit m}}{H^{b}_{(p,b)(p,b)}} (g^{(t,a)(p,b)})^2
\end{align}
If $a=b$, then the r.h.s of \eqref{H_A^m_cal} equals
\begin{align*}
&\sum_{t>a}(\mathcomposite[3]{\mathrel}{\hspace{-11mm} {\scriptsize\textit m}}{H^{a}_{(a,t)(a,t)}} )^2 (g^{(a,t)(a,t)})^2+
\sum_{t<a}(\mathcomposite[3]{\mathrel}{\hspace{-11mm} {\scriptsize\textit m}}{H^{a}_{(t,a)(t,a)}} )^2 (g^{(t,a)(t,a)})^2\\
&=\sum_{t>a}(2(\lambda_t-\lambda_a))^2\Bigl(\frac{\lambda_a \lambda_t}{(\lambda_a-\lambda_t)^2}\Bigr)^2
+\sum_{t<a}(2(\lambda_t-\lambda_a))^2\Bigl(\frac{\lambda_a \lambda_t}{(\lambda_a-\lambda_t)^2}\Bigr)^2\\
&=4\sum_{t \ne a}\frac{\lambda_a^2\lambda_t^2}{(\lambda_a-\lambda_t)^2}.
\end{align*}
If $a\ne b$, then the r.h.s of \eqref{H_A^m_cal} equals
\begin{align*}
&\mathcomposite[3]{\mathrel}{\hspace{-11mm} {\scriptsize\textit m}}{H^{a}_{(a,b)(a,b)}} \:\mathcomposite[1]{\mathrel}{\hspace{-11mm} {\scriptsize\textit m}}{H^{b}_{(a,b)(a,b)}} (g^{(a,b)(a,b)})^2\\
&=4(\lambda_b-\lambda_a)(\lambda_a-\lambda_b)\Bigl( \frac{\lambda_a\lambda_b}{(\lambda_a-\lambda_b)^2}\Bigr)^2\\
&=-\frac{4\lambda_a^2\lambda_b^2}{(\lambda_a-\lambda_b)^2}.
\end{align*}
\subsection{Proof of Proposition \ref{info_loss_result}}
The term of the order $n$ in \eqref{information_loss} vanishes since $g_{a(s,t)}$ equals zero for $1\leq a \leq p,\ 1\leq s<t \leq p$. We consider the term of order $O(1)$. Since $\stackrel{e}{H}_{ac(s,t)}$ also vanishes for $1\leq a, c \leq p$, $1\leq s<t \leq p$, we only have to consider the term
$$
(1/2)\sum_{s<t,u<v,o<p,q<r}\stackrel{m}{H}_{(s,t)(u,v)a}\;\stackrel{m}{H}_{(o,p)(q,r)b}\;g^{(s,t)(o,p)}\;g^{(u,v)(q,r)}.
$$
Because of \eqref{embedding_curvature_A}, the above term equals
\begin{align}
\label{info_loss_cal}
&2^{-1}\sum_{t>a, p>b} \stackrel{m}{H}_{(a,t)(a,t)a}\stackrel{m}{H}_{(b,p)(b,p)b}(g^{(a,t)(b,p)})^2 \nonumber\\
&+2^{-1}\sum_{t>a, p<b} \stackrel{m}{H}_{(a,t)(a,t)a}\stackrel{m}{H}_{(p,b)(p,b)b}(g^{(a,t)(p,b)})^2 \nonumber\\
&+2^{-1}\sum_{t<a, p>b} \stackrel{m}{H}_{(t,a)(t,a)a}\stackrel{m}{H}_{(b,p)(b,p)b}(g^{(t,a)(b,p)})^2 \nonumber\\
&+2^{-1}\sum_{t<a, p<b} \stackrel{m}{H}_{(t,a)(t,a)a}\stackrel{m}{H}_{(p,b)(p,b)b}(g^{(t,a)(p,b)})^2
\end{align}
If $a=b$, then \eqref{info_loss_cal} equals 
\begin{align*}
&2^{-1}\sum_{t>a} (\stackrel{m}{H}_{(a,t)(a,t)a})^2(g^{(a,t)(a,t)})^2
+2^{-1}\sum_{t<a} (\stackrel{m}{H}_{(t,a)(t,a)a})^2(g^{(t,a)(t,a)})^2\\
&=2^{-1}\Bigl\{\sum_{t>a}(\lambda_a^{-2}(\lambda_t-\lambda_a))^2\Bigl(\frac{\lambda_a \lambda_t}{(\lambda_a-\lambda_t)^2}\Bigr)^2
+\sum_{t<a}(\lambda_a^{-2}(\lambda_t-\lambda_a))^2\Bigl(\frac{\lambda_a \lambda_t}{(\lambda_a-\lambda_t)^2}\Bigr)^2\Bigr\}\\
&=2^{-1}\sum_{t \ne a}\frac{\lambda_t^2}{\lambda_a^2(\lambda_a-\lambda_t)^2}.
\end{align*}
If $a < b$, then \eqref{info_loss_cal} equals 
\begin{align*}
&2^{-1} \stackrel{m}{H}_{(a,b)(a,b)a}\stackrel{m}{H}_{(a,b)(a,b)b}(g^{(a,b)(a,b)})^2\\
&=2^{-1}\lambda_a^{-2}(\lambda_b-\lambda_a)\lambda_b^{-2}(\lambda_a-\lambda_b)\Bigl( \frac{\lambda_a\lambda_b}{(\lambda_a-\lambda_b)^2}\Bigr)^2\\
&=-\frac{1}{2(\lambda_a-\lambda_b)^2}.
\end{align*}
\section*{Acknowledgment}
The author really appreciates Dr. M. Kumon kindly answering his question on a basic fact of the information loss.
He also shows deep gratitude to anonymous referees for their valuable comments which improved the quality of the paper.

\end{document}